\input amstex
\documentstyle{amsppt}
\accentedsymbol\tx{\tilde x}
\accentedsymbol\ty{\tilde y}
\def\const{\operatorname{const}}
\def\sll{\operatorname{sl}}
\pagewidth{12.9cm}
\pageheight{21.3cm}
\nopagenumbers

\topmatter
\title
Effective procedure of point-classification for the 
equation $y''=P+3\,Q\,y'+3\,R\,{y'}^2+S\,{y'}^3$.
\endtitle
\rightheadtext{Effective procedure of point-classification \dots}
\author
R.~A.~Sharipov
\endauthor
\thanks
Paper is written under financial support of European fund INTAS
(project \#93-47, coordinator of project S.~I.~Pinchuk) and
Russian fund for Fundamental Researches (project \#96-01-00127, 
head of project Ya.~T.~Sultanaev). Work is also supported by the
grant of Academy of Sciences of the Republic Bashkortostan
(head of project N.~M.~Asadullin).
\endthanks
\address
Rabochaya str. 5, Ufa 450003, Russia.
\endaddress
\email
root\@bgua.bashkiria.su
\endemail
\abstract\nofrills For the equations of the form $y''=P(x,y)+
3\,Q(x,y)\,y'+3\,R(x,y)\,{y'}^2+S(x,y)\,{y'}^3$ the problem of
equivalence in the class of point transformations is considered.
Effective procedure for determining the class of point 
equivalence for the given equation is suggested. This procedure
is based on explicit formulas for the invariants.
\endabstract
\endtopmatter
\loadbold
\document
\head
1. Introduction.
\endhead
    Let's consider an ordinary differential equation $y''=f(x,y,y')$,
where $f(x,y,y')$ is the third order polynomial in $y'$:
$$
y''=P(x,y)+3\,Q(x,y)\,y'+3\,R(x,y)\,{y'}^2+S(x,y)\,{y'}^3.
\tag1.1
$$
We apply the following point-transformation to it:
$$
\cases
\tilde x=\tilde x(x,y),\\
\tilde y=\tilde y(x,y).
\endcases
\tag1.2
$$
As a result of such transformation we get another equation
which has the same form:
$$
\ty''=\tilde P(\tx,\ty)+3\,\tilde Q(\tx,\ty)\,\ty'+3\,\tilde R(\tx,\ty)
\,{\ty'}^2+\tilde S(\tx,\ty)\,{\ty'}^3.
\tag1.3
$$
Two equations \thetag{1.1} and \thetag{1.3} which are 
bound with the transformation \thetag{1.2} are called {\it 
point-equivalent equations}. The problem of finding criteria
for detecting point-equivalence for two given equations
\thetag{1.1} and \thetag{1.3} is known as {\it the problem of
equivalence}. This problem was studied in numerous papers
(see \cite{1--23}), some of them are classical papers and others 
are modern ones. Results of these papers were summed up in
\cite{24}. In that paper the complete description of 
point-equivalence classes for the equations of the form 
\thetag{1.1} is given.\par
     Some special classifying parameters play the key role in
describing point-equi\-va\-lence classes. Some of them are scalar
invariants for the equation \thetag{1.1}, others are the
components of pseudotensorial fields of various weights.
\definition{Definition 1.1} Pseudotensorial field of the type
$(r,s)$ and weight $m$ is an indexed array $F^{i_1\ldots\,i_r}_{j_1
\ldots\,j_s}$ that under the point transformations \thetag{1.2}
is transformed as follows:
$$
F^{i_1\ldots\,i_r}_{j_1\ldots\,j_s}=
(\det T)^m\sum\Sb p_1\ldots p_r\\ q_1\ldots q_s\endSb
S^{i_1}_{p_1}\ldots\,S^{i_r}_{p_r}\,\,
T^{q_1}_{j_1}\ldots\,T^{q_s}_{j_s}\,\,
\tilde F^{p_1\ldots\,p_r}_{q_1\ldots\,q_s}.
\tag1.4
$$
Here in \thetag{1.4} $T$ and $S$ are Jacobian matrices\footnotemark\ for 
direct and inverse point transformations respectively:
$$
\xalignat 2
&S=\Vmatrix
x_{\sssize 1.0} &x_{\sssize 0.1}\\
\vspace{1ex}
y_{\sssize 1.0} &y_{\sssize 0.1}
\endVmatrix,
&&T=\Vmatrix
\tilde x_{\sssize 1.0} &\tilde x_{\sssize 0.1}\\
\vspace{1ex}
\tilde y_{\sssize 1.0} &\tilde y_{\sssize 0.1}
\endVmatrix.
\tag1.5
\endxalignat
$$
\enddefinition
\noindent First two classifying parameters are determined by
the coefficients of the equation \thetag{1.1} according to
the following formulas:
$$
\aligned
&\aligned
 A=P_{\sssize 0.2}&-2\,Q_{\sssize 1.1}+R_{\sssize 2.0}+
 2\,P\,S_{\sssize 1.0}+S\,P_{\sssize 1.0}-\\
 \vspace{0.5ex}
 &-3\,P\,R_{\sssize 0.1}-3\,R\,P_{\sssize 0.1}-
 3\,Q\,R_{\sssize 1.0}+6\,Q\,Q_{\sssize 0.1},
 \endaligned\\
 \vspace{1ex}
&\aligned
 B=S_{\sssize 2.0}&-2\,R_{\sssize 1.1}+Q_{\sssize 0.2}-
 2\,S\,P_{\sssize 0.1}-P\,S_{\sssize 0.1}+\\
 \vspace{0.5ex}
 &+3\,S\,Q_{\sssize 1.0}+3\,Q\,S_{\sssize 1.0}+
 3\,R\,Q_{\sssize 0.1}-6\,R\,R_{\sssize 1.0}.
 \endaligned
\endaligned
\tag1.6
$$
Parameters $A$ and $B$ are the components of pseudovectorial field
$\alpha$ of the weight two: $\alpha^1=B$ and $\alpha^2=-A$.
The case when both these parameters are zero is known as {\it the
case of maximal degeneration} (see \cite{24}):\footnotetext{By means
of double indices in \thetag{1.5} and in what follows we denote
partial derivatives. Thus for the function $f(x,y)$ by $f_{\sssize p.q}$
we denote the differentiation $p$-times with respect to $x$ and
$q$-times with respect to $y$.}
$$
\xalignat 2
&A=0,&&B=0
\tag1.7
\endxalignat
$$
This case is well known. Each equation for which the conditions 
\thetag{1.7} hold is point-equivalent to the trivial one $y''=0$.
Each such equation has eight-parametric group of point symmetries
isomorphic to $SL(3,\Bbb R)$.\par
     The conditions \thetag{1.7} that determine the case of maximal
degeneration are absolutely effective. In order to check these 
conditions one should only differentiate coefficients of the 
equation \thetag{1.1} and substitute them and their derivatives 
into \thetag{1.6}. But, apart from the case of maximal degeneration,
the complete scheme of point-classification from \cite{24} includes
eight more cases: {\it the case of general position} and seven
cases of {\it intermediate degeneration}. Conditions that determine
these cases are much less effective. The matter is that these
conditions are formulated in special variables, where the 
pseudovectorial field $\alpha$ has unitary components:
$$
\xalignat 2
&\alpha^1=B=1,&&\alpha^2=-A=0.
\tag1.8
\endxalignat
$$
Components of any nonzero pseudovectorial field of the weight
$2$ can be brought to the form \thetag{1.8} in some variables. 
But for to find appropriate variables one should solve some
system of ordinary differential equation. Theoretically this
make no limitations and the procedure of point-classification
was completed in \cite{24}. In practice this step produces
the great deal of inefficiency, since explicitly solvable
system of ordinary differential equations is very rare event. 
Main goal of this paper is to eliminate this inefficient step 
of bringing $\alpha$ into the form \thetag{1.8}. Then the 
point-classification procedure from paper \cite{24} will 
become absolutely effective.\par
     It is worth to note that analogous inefficiency took
place in the theory of hydrodynamical equations which are 
integrable by generalized hodograph method. This inefficiency 
was eliminated in \cite{25} when the integrability condition 
was written in an invariant tensorial form.
\head
2. Pseudoscalar field $F$ and pseudovectorial field $\beta$.
\endhead
     According to the point-classification scheme from \cite{24},
in all cases different from the case of maximal degeneration the
study of the equation \thetag{1.1} starts from the construction
of pseudoscalar field $F$ of the weight $1$. In special 
coordinates, where the conditions \thetag{1.8} hold, this field
is given by the formula
$$
F^5=-P.
\tag2.1
$$
The fact that the formula \thetag{2.1} gives the pseudoscalar
field of the weight $1$ was checked in \cite{24}. Now we are
only to recalculate this field in arbitrary coordinates.
\proclaim{Theorem 2.1} Pseudoscalar field $F$ of the weight $1$
in arbitrary coordinates is given by the formula
$$
\aligned
F^5=A\,B\,A_{\sssize 0.1}&+B\,A\,B_{\sssize 1.0}-
A^2\,B_{\sssize 0.1}-B^2\,A_{\sssize 1.0}-\\
&-P\,B^3+3\,Q\,A\,B^2-3\,R\,A^2\,B+S\,A^3,
\endaligned
\tag2.2
$$
where parameters $A$ and $B$ should be calculated by the
formulas \thetag{1.6}.
\endproclaim
    It is easy to check that when we substitute \thetag{1.8}
into \thetag{2.2}, this formula reduces to \thetag{2.1}.
Therefore in order to prove the theorem~2.1 we have only to
prove that the formula \thetag{2.2} determines the pseudoscalar
field of the weight $1$. We shall not do it here, since this
was done in \cite{26}. The formula \thetag{2.2} itself was
also derived in \cite{26}.\par
     Together with $F$, in \cite{24} and \cite{26} the pseudovectorial
field $\beta$ of the weight $4$ was defined. In special coordinates 
it has the following components:
$$
\xalignat 2
&\beta^1=G=3\,Q,&&\beta^2=H=-3\,P.
\tag2.3
\endxalignat
$$
In arbitrary coordinates formulas \thetag{2.3} are recalculated
into the following ones:
$$
\aligned
G&=-B\,B_{\sssize 1.0}-3\,A\,B_{\sssize 0.1}+4\,B\,
A_{\sssize 0.1}+3\,S\,A^2-6\,R\,B\,A+3\,Q\,B^2,\\
\vspace{1ex}
H&=-A\,A_{\sssize 0.1}-3\,B\,A_{\sssize 1.0}+4\,A\,
B_{\sssize 1.0}-3\,P\,B^2+6\,Q\,A\,B-3\,R\,A^2.
\endaligned
\tag2.4
$$
For raising and lowering indices we shall use the skew-symmetric 
matrix
$$
d_{ij}=d^{ij}=
\Vmatrix\format \r&\quad\l\\ 0 & 1\\-1 & 0\endVmatrix.
\tag2.5
$$
Components of the matrix $d_{ij}$ form twice-covariant pseudotensorial
field of the weight $-1$, the same values denoted by $d^{ij}$ form
twice-contravariant field of the weight $1$. Pseudotensorial fields
\thetag{2.5} reveal the relationship between $\alpha$, $\beta$ and $F$:
$$
3\,F^5=\sum^2_{i=1}\alpha_i\,\beta^i=\sum^2_{i=1}\sum^2_{j=1}
d_{ij}\,\alpha^i\,\beta^j.
\tag2.6
$$
    All above statements concerning the fields $\alpha$, $\beta$ and $F$,
as well as the explicit formulas for them, are derived on the base of
transformation rules for the coefficients of the equation \thetag{1.1}.
In order to write these rules in a brief form let's construct the 
following three-index array:
$$
\xalignat 2
&\theta_{111}=P,
&&\theta_{112}=\theta_{121}=\theta_{211}=Q,
\hskip -2em\\
\vspace{-1.7ex}
&&&\tag2.7\\
\vspace{-1.7ex}
&\theta_{122}=\theta_{212}=\theta_{221}=R,
&&\theta_{222}=S.
\hskip -2em
\endxalignat
$$
Then let's raise one of the indices by means of the matrix $d^{ij}$:
$$
\theta^k_{ij}=\sum^2_{r=1}d^{kr}\,\theta_{rij}.
\tag2.8
$$
Under the transformation \thetag{1.2} the components of the array 
$\theta^k_{ij}$ are transformed as
$$
\theta^k_{ij}=\sum^2_{m=1}\sum^2_{p=1}\sum^2_{q=1}
S^k_m\,T^p_i\,T^q_j\,\tilde\theta^m_{pq}
+\sum^2_{m=1}S^k_m\,\frac{\partial T^m_i}{\partial x^j}-
\frac{\tilde\sigma_i\,\delta^k_j+
\tilde\sigma_j\,\delta^k_i}{3}.
\tag2.9
$$
Here $x^1=x$, $x^2=y$, $\tilde x^1=\tilde x$, $\tilde x^2=\tilde y$,
and, moreover, the following notations are made:
$$
\xalignat 3
&\qquad\tilde\sigma_i=\frac{\partial\ln\det T}{\partial x^i},
&&\delta^k_i=
\cases
1&\text{for\ }i=k,\\
0&\text{for\ }i\neq k.
\endcases
\hskip -2em
\endxalignat
$$
Transformation rules \thetag{2.9} for the components of the array
\thetag{2.8} are quite similar to that for the components of an 
affine connection (see \cite{27}). The difference consists only in 
the fraction with the number $3$ in denominator.
\head
3. The case of general position.
\endhead
     According to the scheme of point-classification from \cite{24},
the case of general position is determined by the condition $F\neq 
0$. This condition is equivalent to the non-collinearity of
pseudovectorial fields $\alpha$ and $\beta$. This fact is easily
derived from the formula \thetag{2.6}. Moreover, the condition 
$F\neq 0$ makes possible to define the following logarithmic 
derivatives:
$$
\varphi_i=-\frac{\partial\ln F}{\partial x^i}.
\tag3.1
$$
Under the point transformations \thetag{1.2} the quantities
\thetag{3.1} are transformed as
$$
\varphi_i=\sum^2_{j=1}T^j_i\,\tilde\varphi_j-\tilde\sigma_i.
\tag3.2
$$
Relying on \thetag{3.2}, we can use the quantities $\varphi_i$
and the quantities \thetag{2.8} in order to construct the 
components of an affine connection:
$$
\varGamma^k_{ij}=\theta^k_{ij}-\frac{\varphi_i\,\delta^k_j+
\varphi_j\,\delta^k_i}{3}.
\tag3.3
$$
In addition to \thetag{3.3}, we can define the pair of vectorial
fields
$$
\xalignat 2
&\bold X=F^{-2}\,\alpha,
&&\bold Y=F^{-4}\,\beta.
\tag3.4
\endxalignat
$$
The same condition $F\neq 0$ warrants the non-collinearity
of the fields $\bold X$ and $\bold Y$ from \thetag{3.4}. Thus
we have the moving frame in the plane of variables $x$ and $y$.
It is formed by vectors $\bold X$ and $\bold Y$ at each point.
Connection components \thetag{3.3} define the covariant
differentiation of vector fields. Now let's calculate the
components of the connection \thetag{3.3} referred to the
frame of the fields $\bold X$ and $\bold Y$. They are defined
as the coefficients in the following expansions:
$$
\xalignat 2
&\nabla_{\bold X}\bold X=\Gamma^1_{11}\,\bold X+
\Gamma^2_{11}\,\bold Y,
&
&\nabla_{\bold X}\bold Y=\Gamma^1_{12}\,\bold X+
\Gamma^2_{12}\,\bold Y,\\
\vspace{-1.7ex}
&&&\tag3.5\\
\vspace{-1.7ex}
&\nabla_{\bold Y}\bold X=\Gamma^1_{21}\,\bold X+
\Gamma^2_{21}\,\bold Y,
&
&\nabla_{\bold Y}\bold Y=\Gamma^1_{22}\,\bold X+
\Gamma^2_{22}\,\bold Y.
\endxalignat
$$
The quantities $\Gamma^k_{ij}$ in the expansions \thetag{3.5}
are the scalar invariants of the equation \thetag{1.1}. In
paper \cite{24} they were denoted as follows:
$$
\xalignat 4
&\Gamma^1_{11}=I_1,&&\Gamma^2_{11}=I_2,
&&\Gamma^1_{12}=I_3,&&\Gamma^2_{12}=I_4,\\
&\Gamma^1_{21}=I_5,&&\Gamma^2_{21}=I_6,
&&\Gamma^1_{22}=I_7,&&\Gamma^2_{22}=I_8.
\endxalignat
$$
By differentiating these invariants along vector field 
$\bold X$ we get eight more invariants
$$
\xalignat 4
&\bold XI_1=I_9,&&\bold XI_2=I_{10},&&\bold XI_3=I_{11},
&&\bold XI_4=I_{12},\\
&\bold XI_5=I_{13},&&\bold XI_6=I_{14},&&\bold XI_7=I_{15},&&
\bold XI_8=I_{16}.
\endxalignat
$$
The differentiation of the eight initial invariants along
another vector field increases the number of invariants up
to $24$:
$$
\xalignat 4
&\bold YI_1=I_{17},&&\bold YI_2=I_{18},&&\bold YI_3=I_{19},
&&\bold YI_4=I_{20},\\
&\bold YI_5=I_{21},&&\bold YI_6=I_{22},&&\bold YI_7=I_{23},&&
\bold YI_8=I_{24}.
\endxalignat
$$
Repeating this procedure of differentiation along $\bold X$ and
$\bold Y$, we can construct indefinite sequence of invariants
by adding $16$ ones in each step. According to the results of
\cite{24} and \cite{26}, the properties of this sequence of
invariants divide the case of general position into three
subcases:
\roster
\item in the infinite sequence of invariants $I_k(x,y)$ one can
      find two functionally independent ones;
\item invariants $I_k(x,y)$ are functionally dependent, but not
      all of them are constants;
\item all invariants in the sequence $I_k(x,y)$ are constants.
\endroster
In the first case the group of point symmetries of the equation
\thetag{1.1} is trivial, in the second case it is one-dimensional,
and in the third case it is two-dimensional (see theorems~5.1, 
5.2, and 5.3 in \cite{24}).\par
     In order to tell which of these three cases takes place for
the given particular equation of the form \thetag{1.1} we should
calculate invariants $I_k$ explicitly. Note that we can do this
effectively, since the components of the vector fields $\bold X$ 
and $\bold Y$ and the components of connection \thetag{3.3} are 
now calculated in arbitrary coordinates. We shall give the
explicit formulas for eight first invariants in the sequence.
These formulas were obtained in \cite{26}:
$$
\align
&\hskip -1em\aligned
 I_3&=\frac{B\,(H\,G_{\sssize 1.0}-G\,H_{\sssize 1.0})}{3\,F^9}
 -\frac{A\,(H\,G_{\sssize 0.1}-G\,H_{\sssize 0.1})}{3\,F^9}+
 \frac{H\,F_{\sssize 0.1}+G\,F_{\sssize 1.0}}{3\,F^5}+\\
 \vspace{1ex}
 &+\frac{B\,G^2\,P}{3\,F^9}-\frac{(A\,G^2-2\,H\,B\,G)\,Q}{3\,F^9}+
 \frac{(B\,H^2-2\,H\,A\,G)\,R}{3\,F^9}-\frac{A\,H^2\,S}{3\,F^9},
 \endaligned\hskip -4em
\tag3.6\\
\vspace{3ex}
&\hskip -1em\aligned
 I_6&=\frac{A\,(G\,A_{\sssize 0.1}+H\,B_{\sssize 0.1})}{12\,F^7}
 -\frac{B\,(G\,A_{\sssize 1.0}+H\,B_{\sssize 1.0})}{12\,F^7}
 -\frac{4\,(A\,F_{\sssize 0.1}-B\,F_{\sssize 1.0})}{12\,F^3}-\\
 \vspace{1ex}
 &-\frac{G\,B^2\,P}{12\,F^7}-\frac{(H\,B^2-2\,G\,B\,A)\,Q}{12\,F^7}-
 \frac{(G\,A^2-2\,H\,B\,A)\,R}{12\,F^7}+\frac{H\,A^2\,S}{12\,F^7}.
 \endaligned\hskip -4em
\tag3.7
\endalign
$$
Seventh invariant $I_7$ is given by the formula
$$
\aligned
 I_7=&\frac{G\,H\,G_{\sssize 1.0}-
 G^2\,H_{\sssize 1.0}+H^2\,G_{\sssize 0.1}-
 H\,G\,H_{\sssize 0.1}}{3\,F^{11}}+\\
 \vspace{1ex}
 &\hskip 5em+\frac{G^3\,P+3\,G^2\,H\,\,Q+3\,G\,H^2\,R+H^3\,S}
 {3\,F^{11}}.
 \endaligned
\tag3.8
$$
Formula for the eighth invariant $I_8$ is similar to \thetag{3.6}
and \thetag{3.7}:
$$
\aligned
I_8&=\frac{G\,(A\,G_{\sssize 1.0}+B\,H_{\sssize 1.0})}{3\,F^9}
 +\frac{H\,(A\,G_{\sssize 0.1}+B\,H_{\sssize 0.1})}{3\,F^9}
 -\frac{10\,(H\,F_{\sssize 0.1}+G\,F_{\sssize 1.0})}{3\,F^5}-
 \vspace{1ex}
 &-\frac{B\,G^2\,P}{3\,F^9}+\frac{(A\,G^2-2\,H\,B\,G)\,Q}{3\,F^9}-
 \frac{(B\,H^2-2\,H\,A\,G)\,R}{3\,F^9}+\frac{A\,H^2\,S}{3\,F^9}.
 \endaligned\hskip -1em
\tag3.9
$$
Invariants $I_1$, $I_2$, $I_4$, and $I_5$ do not require separate
calculations. They are expressed through that ones, which are
already calculated:
$$
\xalignat 4
&\quad I_1=-4\,I_6,&&I_2=\frac{1}{3},&&I_4=4\,I_6,&&I_5=-I_8.
\hskip -2em
\tag3.10
\endxalignat
$$
Note that the formula \thetag{3.7} for $I_6$ can be simplified
substantially:
$$
I_6=\frac{A_{\sssize 0.1}-B_{\sssize 1.0}}{3\,F^2}-
\frac{A\,F_{\sssize 0.1}-B\,F_{\sssize 1.0}}{3\,F^3}.
\tag3.11
$$
Due to \thetag{3.10} the formula \thetag{3.11} simplifies the
calculation for two more invariants $I_1$ and $I_4$. As for the
simplification of the formulas \thetag{3.6}, \thetag{3.8}, and 
\thetag{3.9} we shall not undertake special efforts for this now.
\head
4. Cases of intermediate degeneration.
\endhead
    Apart from the case of general position and the case of maximal
degeneration, the scheme of point-classification from \cite{24}
includes seven cases of {\it intermediate degeneration} when
parameters \thetag{1.6} does not vanish simultaneously. All them
are characterized by the condition $F=0$. Due to the formula
\thetag{2.2} this condition is written as follows:
$$
\aligned
A\,B\,A_{\sssize 0.1}&+B\,A\,B_{\sssize 1.0}-
A^2\,B_{\sssize 0.1}-B^2\,A_{\sssize 1.0}-\\
&-P\,B^3+3\,Q\,A\,B^2-3\,R\,A^2\,B+S\,A^3=0.
\endaligned
\tag4.1
$$
The condition \thetag{4.1} is equivalent to the collinearity
of pseudovectorial fields $\alpha$ and $\beta$. The field
$\alpha$ is nonzero (since otherwise we would have the case of
maximal degeneration). Therefore the condition of collinearity 
$\alpha\parallel\beta$ can be written as
$$
\beta=3\,N\,\alpha.
\tag4.2
$$
The field $\beta$ has the weight $4$, while the field $\alpha$
is of the weight $2$. Therefore the coefficient of proportionality
$N$ in \thetag{4.2} is the pseudoscalar field of the weight $2$. 
Pseudoscalar field $N$ defined by \thetag{4.1} can be evaluated 
by any one of the following two formulas:
$$
\xalignat 2
&N=\frac{G}{3\,B},&&N=-\frac{H}{3\,A},
\tag4.3
\endxalignat
$$
Here $A$, $B$, $G$ and $H$ are calculated by \thetag{1.6} and 
\thetag{2.4}. In case of vanishing either $A$ or $B$ one of the 
formulas \thetag{4.3} gives uncertainty $0/0$, but the other
remains true for to calculate $N$.\par
     In special coordinates, where the conditions \thetag{1.8}
hold, the field $N$ is calculated by the first formula \thetag{4.3}.
Here we have 
$$
\xalignat 2
&P=0,&& N=Q.
\tag4.4
\endxalignat
$$
The relationships \thetag{4.4} are in complete agreement with the 
results of \cite{24}. In that paper, besides $N=Q$, the following 
two pseudoscalar fields were defined:
$$
\xalignat 2
&\Omega=R_{\sssize 1.0}-2\,Q_{\sssize 0.1},
&&M=Q_{\sssize 1.0}-\frac{12}{5}\,Q^2.
\tag4.5
\endxalignat
$$
The field $\Omega$ is of the weight $1$, and $M$ is of the weight $4$.
Pseudovectorial field $\gamma$ has the weight $3$:
$$
\xalignat 2
&\gamma^1=-2\,R_{\sssize 1.0}+3\,Q_{\sssize 0.1}+
\frac{6}{5}\,Q\,R,
&&\gamma^2=M.
\tag4.6
\endxalignat
$$
Due to $F=0$ the quantities \thetag{3.1} in all cases of intermediate
degeneration are not defined. Instead of them in \cite{24} other
two quantities $\varphi_1$ and $\varphi_2$ were introduced: 
$$
\xalignat 2
&\varphi_1=-\frac{6}{5}\,Q,
&&\varphi_2=-\frac{3}{5}\,R.
\tag4.7
\endxalignat
$$
They obey the same transformation rule \thetag{3.2} as the quantities 
\thetag{3.1}. Formulas \thetag{4.5}, \thetag{4.6} and \thetag{4.7},
taken from \cite{24}, are written in special coordinates. In order
to make effective we should recalculate them in arbitrary
coordinates. Now we can't use the results of \cite{26}, since the
cases of intermediate degeneration are not considered there.
\proclaim{Theorem 4.1} If $B\neq 0$, then in any case of
intermediate degeneration the parameters $\varphi_i$ are defined 
by the formulas
$$
\aligned
&\varphi_1=-3\,A\,\frac{A\,S-B_{\sssize 0.1}}{5\,B^2}
-3\,\frac{A_{\sssize 0.1}+B_{\sssize 1.0}-3\,A\,R}{5\,B}
-\frac{6}{5}\,Q,\\
\vspace{2ex}
&\varphi_2=3\,\frac{A\,S-B_{\sssize 0.1}}{5\,B}
-\frac{3}{5}\,R,
\endaligned
\tag4.8
$$
which hold for arbitrary curvilinear coordinates $x$ and $y$ 
on the plane.
\endproclaim
    It is easy to note that in special coordinates, where the
conditions \thetag{1.8} hold, the relationships \thetag{4.8} 
are reduced to the form \thetag{4.7}.
\demo{Proof} Let $\tx$ and $\ty$ be special coordinates, for
which the conditions \thetag{1.8} are fulfilled. This can be
expressed as follows:
$$
\xalignat 2
&\tilde\alpha^1=\tilde B=1,&&\tilde\alpha^2=-\tilde A=0.
\endxalignat
$$
Then for to calculate the parameters $A$ and $B$ in arbitrary
(nonspecial) coordinates $x$ and $y$ one can use the pseudovectorial
rule of transformation for the components of the field $\alpha$.
Using definition~1.1, we get
$$
\aligned
B&=\alpha^1=(\det T)^2(S^1_1\,\tilde\alpha^1+S^1_2\,\tilde\alpha^2)=
\ty_{\sssize 0.1}^2\,\tx_{\sssize 1.0}-\ty_{\sssize 1.0}\,
\ty_{\sssize 0.1}\,\tx_{\sssize 0.1},\\
\vspace{1ex}
-A&=\alpha^2=(\det T)^2(S^2_1\,\tilde\alpha^1+S^2_2\,\tilde\alpha^2)=
\ty_{\sssize 1.0}^2\,\tx_{\sssize 0.1}-\ty_{\sssize 1.0}\,
\ty_{\sssize 0.1}\,\tx_{\sssize 1.0}.
\endaligned
\hskip-2em
\tag4.9
$$
From $B\neq 0$ and from \thetag{4.9} we conclude that $\ty_{\sssize 0.1}
\neq 0$. Therefore we can resolve the equations \thetag{4.9} with
respect to the derivatives $\tx_{\sssize 1.0}$ and $\ty_{\sssize 1.0}$:
$$
\xalignat 2
&\tx_{\sssize 1.0}=\frac{A}{B}\,\ty_{\sssize 0.1},
&&\ty_{\sssize 1.0}=\frac{B}{\ty_{\sssize 0.1}^2}+\frac{A}{B}
\,\tx_{\sssize 0.1}.
\hskip -2em
\tag4.10
\endxalignat
$$
By differentiating \thetag{4.10} with respect to $x$ and $y$
we get formulas that express $\tx_{\sssize 2.0}$, $\tx_{\sssize 1.1}$, 
$\ty_{\sssize 2.0}$, and $\ty_{\sssize 1.1}$ through $\tx_{\sssize 0.1}$ 
and $\ty_{\sssize 0.1}$. Substituting the obtained expressions for
second order derivatives into the transformation rules for the
coefficients of the equation \thetag{1.1}, we derive the following
relationship for $S$:
$$
S=\frac{\ty_{\sssize 0.1}^2\,\tx_{\sssize 0.2}}{B}
-\frac{\ty_{\sssize 0.1}\,\ty_{\sssize 0.2}\,\tx_{\sssize 0, 1}}{B}
+3\,\frac{\ty_{\sssize 0.1}^2\,\tilde Q\,\tx_{\sssize 0.1}^2}{B}
+\frac{\ty_{\sssize 0.1}^4\,\tilde S}{B}
+3\,\frac{\ty_{\sssize 0.1}^3\,\tilde R\,\tx_{\sssize 0.1}}{B}.
\hskip -2em
\tag4.11
$$
Now let's use the relationship \thetag{4.11} for to express
the derivative $\tx_{\sssize 0.2}$ through $\tx_{\sssize 0.1}$
and $\ty_{\sssize 0.1}$. Upon doing this, for the coefficient
$R$ we obtain
$$
R=-\frac{5}{3}\,\frac{\ty_{\sssize 0.2}}{\ty_{\sssize 0.1}}
+\frac{2}{3}\,\frac{B_{\sssize 0.1}}{B}+\frac{A\,S}{B}
+\ty_{\sssize 0.1}\,\tilde R+2\,\tilde Q\,\tx_{\sssize 0.1}.
\hskip -2em
\tag4.12
$$
The relationship \thetag{4.12} can be used to express 
$\ty_{\sssize 0.2}$ through $\tx_{\sssize 0.1}$ and 
$\ty_{\sssize 0.1}$. Then formula for $Q$ can be brought
to the form
$$
Q=-\frac{4}{3}\,\frac{A_{\sssize 0.1}}{B}+\frac{1}{3}
\,\frac{B_{\sssize 1.0}}{B}+\frac{B\,\tilde Q}{\ty_{\sssize 0.1}^2}
+\frac{A\,B_{\sssize 0.1}}{B^2}-\frac{A^2\,S}{B^2}+2\,\frac{A\,R}{B}.
$$
This formula is used to express $\tilde Q$ through the derivative
$\ty_{\sssize 0.1}$. Then the formula for $P$ can be written as
$$
P=-\frac{A_{\sssize 1.0}}{B}+\frac{A^3\,S}{B^3}
+\frac{A\,A_{\sssize 0.1}}{B^2}+\frac{A\,B_{\sssize 1,0}}
{B^2}-\frac{A^2\,B_{\sssize 0.1}}{B^3}-3\,\frac{A^2\,R}
{B^2}+3\,\frac{A\,Q}{B}.
$$
This formula for $P$ add nothing new, since it can be derived
directly from the condition \thetag{4.1} when $B\neq 0$.
However, the expressions for the second order derivatives
$\tx_{\sssize 2.0}$, $\tx_{\sssize 1.1}$, $\tx_{\sssize 0.2}$, 
$\ty_{\sssize 2.0}$, $\ty_{\sssize 1.1}$, and $\ty_{\sssize 0.2}$ 
are enough to derive the formulas for $\varphi_1$ and 
$\varphi_2$.\par
     For the further calculations we shall use the values of
$\varphi_1$ and $\varphi_2$ in special coordinates. Here we have
$$
\xalignat 2
&\tilde\varphi_1=-\frac{6}{5}\,\tilde Q,
&&\tilde\varphi_2=-\frac{3}{5}\,\tilde R.
\endxalignat
$$
Now let's substitute these values into the transformation rules
\thetag{3.2} for them. We can write \thetag{3.2} as follows:
$$
\varphi_i=\sum^2_{j=1}T^j_i\,\tilde\varphi_j-
\frac{\partial\ln\det T}{\partial x^i}.
\tag4.13
$$
Here $T$ is the transition matrix defined in \thetag{1.5}.
Substituting \thetag{4.10} and all other analogous formulas
for all other second order derivatives $\tx_{\sssize 2.0}$, 
$\tx_{\sssize 1.1}$, $\tx_{\sssize 0.2}$, $\ty_{\sssize 2.0}$, 
$\ty_{\sssize 1.1}$, and $\ty_{\sssize 0.2}$ into \thetag{4.13},
we get the formulas \thetag{4.8}. All occurrences of the 
derivatives $\tx_{\sssize 0.1}$ and $\ty_{\sssize 0.1}$, and
all occurrences of $\tilde S$, $\tilde R$, and $\tilde P$ do
cancel each other during this substitution. Theorem~4.1 is 
proved.\qed\enddemo
    When $B=0$, formulas \thetag{4.8} do not hold. For this
case we have another theorem which can be proved in a similar 
way.
\proclaim{Theorem 4.2} If $A\neq 0$, then in any case of
intermediate degeneration the parameters $\varphi_i$ are defined 
by the formulas
$$
\aligned
&\varphi_1=
-3\,\frac{B\,P+A_{\sssize 1.0}}{5\,A}
+\frac{3}{5}\,Q,\\
\vspace{2ex}
&\varphi_2=
3\,B\,\frac{B\,P+A_{\sssize 1.0}}{5\,A^2}
-3\,\frac{B_{\sssize 1.0}+A_{\sssize 0.1}+3\,B\,Q}{5\,A}
+\frac{6}{5}\,R,
\endaligned
\tag4.14
$$
which hold for arbitrary curvilinear coordinates $x$ and $y$ 
on the plane.
\endproclaim
     For the pseudoscalar field $\Omega$ we need not prove
theorems like theorem~4.1 or theorem~4.2. In paper \cite{24}
one can find invariant definition of this field. According
to \cite{24}, first, we define tensorial field $\omega$
with the following components:
$$
\omega_{ij}=\frac{\partial\varphi_i}{\partial x^j}-
\frac{\partial\varphi_j}{\partial x^i}.
\tag4.15
$$
Then we contract the field \thetag{4.15} with unit skew-symmetric 
field $d$ from \thetag{2.5}:
$$
\Omega=\frac{5}{6}\sum^2_{i=1}\sum^2_{j=1}\omega_{ij}\,
d^{ij}=\frac{5}{3}\left(\frac{\partial\varphi_1}{\partial y}-
\frac{\partial\varphi_2}{\partial x}\right).
\hskip-2em
\tag4.16
$$
Substituting \thetag{4.8} into \thetag{4.16} we get an
explicit formula for $\Omega$ in the case when $B\neq 0$:
$$
\aligned
\Omega&=\frac{2\,A\,B_{\sssize 0.1}(A\,S-B_{\sssize 0.1})}{B^3}
+\frac{(2\,A_{\sssize 0.1}-3\,A\,R)\,B_{\sssize 0.1}}{B^2}+\\
\vspace{1ex}
&+\frac{(B_{\sssize 1.0}-2\,A_{\sssize 0.1})\,A\,S}{B^2}
+\frac{A\,B_{\sssize 0.2}-A^2\,S_{\sssize 0.1}}{B^2}-
\frac{A_{\sssize 0.2}}{B}+\\
\vspace{1ex}
&+\frac{3\,A_{\sssize 0.1}\,R+3\,A\,R_{\sssize 0.1}
-A_{\sssize 1.0}\,S-A\,S_{\sssize 1.0}}{B}
+R_{\sssize 1.0}-2\,Q_{\sssize 0.1}.
\endaligned
\tag4.17
$$
Substituting \thetag{4.8} into \thetag{4.16}, we get another 
formula for $\Omega$, which differs from \thetag{4.17} by the 
following mirror transformation:
$$
\xalignat 2
&x\to y,&&y\to x,\\
&P\to-S,&&S\to -P,\tag4.18\\
&Q\to-R,&&R\to -Q.
\endxalignat
$$
Transformation \thetag{4.18} can be extended to the quantities
$A$, $B$, $G$, $H$, $F$, and $N$ too:
$$
\xalignat 2
&A\to -B,&&B\to -A,\\
&G\to H,&&H\to G,\tag4.19\\
&F\to -F,&&N\to N.
\endxalignat
$$
This can be easily seen from \thetag{2.2}, \thetag{2.4}, and 
\thetag{4.3}. Applying mirror transformations \thetag{4.18} 
and \thetag{4.19} to $\varphi_1$, $\varphi_2$, and $\Omega$,
we get
$$
\xalignat 3
&\varphi_1\to\varphi_2,&&\varphi_2\to\varphi_1,
&&\Omega\to-\Omega.
\hskip -2em
\tag4.20
\endxalignat
$$
On the base of the last transformation in \thetag{4.20} we can
write one more formula for $\Omega$, which holds for $A\neq 0$:
$$
\aligned
\Omega&=\frac{2\,B\,A_{\sssize 1.0}(B\,P+A_{\sssize 1.0})}{A^3}
-\frac{(2\,B_{\sssize 1.0}+3\,B\,Q)\,A_{\sssize 1.0}}{A^2}+\\
\vspace{1ex}
&+\frac{(A_{\sssize 0.1}-2\,B_{\sssize 1.0})\,B\,P}{A^2}
-\frac{B\,A_{\sssize 2.0}+B^2\,P_{\sssize 1.0}}{A^2}+
\frac{B_{\sssize 2.0}}{A}+\\
\vspace{1ex}
&+\frac{3\,B_{\sssize 1.0}\,Q+3\,B\,Q_{\sssize 1.0}
-B_{\sssize 0.1}\,P-B\,P_{\sssize 0.1}}{A}
+Q_{\sssize 0.1}-2\,R_{\sssize 1.0}.
\endaligned
\tag4.21
$$
Formula \thetag{4.21} is mirror symmetric with respect to the 
formula \thetag{4.17} in the sense of the above mirror 
transformations.\par
      Formulas \thetag{4.19} and \thetag{4.21} implement the
effectivization of the formula \thetag{4.5} for $\Omega$. As for
the formulas \thetag{4.5} and \thetag{4.6}, their effectivization
require more efforts. Let's start with the affine
connection given by
$$
\varGamma^k_{ij}=\theta^k_{ij}-\frac{\varphi_i\,\delta^k_j+
\varphi_j\,\delta^k_i}{3}.
\tag4.22
$$
Here the array $\theta^k_{ij}$ is defined by the coefficients of the
equation \thetag{1.1} according to the formulas \thetag{2.7} and 
\thetag{2.8}, while $\varphi_i$ are derived either by \thetag{4.8}
or \thetag{4.14}. We shall use connection components \thetag{4.22} and
the quantities $\varphi_i$ for covariant differentiation of 
pseudotensorial fields:
$$
\aligned
\nabla_kF^{i_1\ldots\,i_r}_{j_1\ldots\,j_s}&=
\frac{\partial F^{i_1\ldots\,i_r}_{j_1\ldots\,j_s}}{\partial x^k}+
\sum^r_{n=1}\sum^2_{v_n=1}\varGamma^{i_n}_{k\,v_n}\,
F^{i_1\ldots\,v_n\ldots\,i_r}_{j_1\ldots\,j_s}-\\
&-\sum^s_{n=1}\sum^2_{w_n=1}\varGamma^{w_n}_{k\,j_n}\,
F^{i_1\ldots\,i_r}_{j_1\ldots\,w_n\ldots\,j_s}+
m\,\varphi_k\,F^{i_1\ldots\,i_r}_{j_1\ldots\,j_s}.
\endaligned
\tag4.23
$$
This formula \thetag{4.23} is the definition of the covariant
derivative $\nabla_k$ for pseudotensorial field $F$ of the
type $(r,s)$ and weight $m$ (for more details see \cite{28}).
As a result of applying \thetag{4.23} we get pseudotensorial
field $\nabla F$ of the type $(r,s+1)$ and weight $m$.\par
     Let's apply the operation of covariant differentiation to the
pseudoscalar field $N$ of the weight $2$. This gives us
pseudocovectorial field $\nabla N$ of the weight $2$. Here are the
components of the field $\nabla N$:
$$
\xalignat 2
&\quad\nabla_1 N=N_{\sssize 1.0}+2\,\varphi_1\,N,
&&\nabla_2 N=N_{\sssize 0.1}+2\,\varphi_2\,N.
\hskip-2em
\tag4.24
\endxalignat
$$
Let's denote by $\xi$ the pseudovectorial field obtained from
$\nabla N$ by raising index by virtue of the matrix $d^{ij}$ 
from \thetag{2.5}:
$$
\xi^i=\sum^2_{j=1}d^{ij}\,\nabla_jN.
\tag4.25
$$
This field $\xi$ has the weight $3$, which coincides with the 
weight of the field $\gamma$. By direct calculations from 
\thetag{4.24} and \thetag{4.25} in special coordinates we get
$$
\xalignat 2
&\quad\xi^1=Q_{\sssize 0.1}-\frac{6}{5}\,R\,Q,
&&\xi^2=-Q_{\sssize 1.0}+\frac{12}{5}\,Q^2=-M.
\hskip-2em
\tag4.26
\endxalignat
$$
Comparing \thetag{4.26} with \thetag{4.5}, \thetag{4.6}, and
\thetag{1.8}, we derive the following relationships between
fields $\alpha$, $\gamma$, $\xi$, $M$, and $\Omega$:
$$
\xalignat 2
&\quad\gamma=-2\,\Omega\,\alpha-\xi,
&&M=-\sum^2_{i=1}\alpha_i\,\xi^i.
\hskip-2em
\tag4.27
\endxalignat
$$
These relationships are written in terms of natural operations of
sum, tensor-product, and contraction for pseudotensorial fields.
Therefore, once they are established in special coordinates, they
remain true in arbitrary coordinates too. For the field $M$ from 
\thetag{4.27} we obtain
$$
\gather
\aligned
M=-\frac{12\,A\,N\,(A\,S-B_{\sssize 0.1})}{5\,B}
&-A\,N_{\sssize 0.1}+\frac{24}{5}\,A\,N\,R-\\
\vspace{1ex}
-\frac{6}{5}\,N\,A_{\sssize 0.1}-\frac{6}{5}\,N\,
&B_{\sssize 1.0}+B\,N_{\sssize 1.0}-\frac{12}{5}\,B\,N\,Q,
\endaligned
\tag4.28\\
\vspace{2ex}
\aligned
M=-\frac{12\,B\,N\,(B\,P+A_{\sssize 1.0})}{5\,A}
&+B\,N_{\sssize 1.0}+\frac{24}{5}\,B\,N\,Q+\\
\vspace{1ex}
+\frac{6}{5}\,N\,B_{\sssize 1.0}+\frac{6}{5}\,N\,
&A_{\sssize 0.1}-A\,N_{\sssize 0.1}-\frac{12}{5}\,A\,N\,R.
\endaligned
\tag4.29
\endgather
$$
Formula \thetag{4.29} is derived from \thetag{4.28} by means
of mirror transformations \thetag{4.18} and \thetag{4.19}. It 
holds for $A\neq 0$, while the initial formula \thetag{4.28}
holds for $B\neq 0$.\par
     For the components of pseudovectorial field $\gamma$ from
\thetag{4.27} we derive the following relationships, which give
the required effectivization for \thetag{4.6} when $B\neq 0$:
$$
\align
&\aligned
\gamma^1=-\frac{6\,N\,(A\,S-B_{\sssize 0.1})}{5\,B}
-N_{\sssize 0.1}+\frac{6}{5}\,N\,R-2\,\Omega\,B,
\endaligned
\tag4.30\\
\vspace{2ex}
&\aligned
\gamma^2=-&\frac{6\,A\,N\,(A\,S-B_{\sssize 0.1})}{5\,B^2}
+\frac{18\,N\,A\,R}{5\,B}-\\
\vspace{1ex}
&-\frac{6\,N\,(A_{\sssize 0.1}+B_{\sssize 1.0})}{5\,B}
+N_{\sssize 1.0}-\frac{12}{5}\,N\,Q+2\,\Omega\,A.
\endaligned
\hskip-2em
\tag4.31
\endalign
$$
Mirror symmetric formulas for $\gamma$, which hold for $A\neq 0$,
have the form
$$
\align
&\aligned
\gamma^1=-&\frac{6\,B\,N\,(B\,P+A_{\sssize 1.0})}{5\,A^2}
+\frac{18\,N\,B\,Q}{5\,A}+\\
\vspace{1ex}
&+\frac{6\,N\,(B_{\sssize 1.0}+A_{\sssize 0.1})}{5\,A}
-N_{\sssize 0.1}-\frac{12}{5}\,N\,R-2\,\Omega\,B,
\endaligned
\hskip-2em
\tag4.32\\
\vspace{2ex}
&\aligned
\gamma^2=-\frac{6\,N\,(B\,P+A_{\sssize 1.0})}{5\,A}
+N_{\sssize 1.0}+\frac{6}{5}\,N\,Q+2\,\Omega\,A.
\endaligned
\tag4.33
\endalign
$$
The mirror transformations themselves for the components of 
$\gamma$ are written as
$$
\xalignat 2
&\gamma^1\to -\gamma^2,&&\gamma^2\to -\gamma^1.
\hskip-2em
\tag4.34
\endxalignat
$$
We can derive \thetag{4.34} by comparing \thetag{4.30} and 
\thetag{4.31} with \thetag{4.32} and \thetag{4.33} and taking 
into account \thetag{4.18} and \thetag{4.19}.
\head
5. First case of intermediate degeneration.
\endhead
     First case of intermediate degeneration is distinguished
from other cases of intermediate degeneration by the condition
$M\neq 0$. This condition $M\neq 0$ is equivalent to the 
condition of non-collinearity $\alpha\nparallel\gamma$ for the
pseudovectorial fields $\alpha$ and $\gamma$. Moreover, from 
$M\neq 0$ one can derive $N\neq 0$. This can be easily seen
either from \thetag{4.28} or from \thetag{4.29}. Therefore
we immediately get two scalar invariants
$$
\xalignat 2
&I_1=\frac{M}{N^2},
&&I_2=\frac{\Omega^2}{N}.
\tag5.1
\endxalignat
$$
Now they can be calculated explicitly in arbitrary coordinates
on the base of the above formulas for $M$, $N$, and $\Omega$. 
In order to define third invariant $I_3$ in \cite{24} the
following expansions were considered:
$$
\xalignat 2
&\nabla_{\alpha}\alpha=\Gamma^1_{11}\,\alpha+
\Gamma^2_{11}\,\gamma,
&&\nabla_{\alpha}\gamma=\Gamma^1_{12}\,\alpha+
\Gamma^2_{12}\,\gamma,\\
\vspace{-1.7ex}
&&&\tag5.2\\
\vspace{-1.7ex}
&\nabla_{\gamma}\alpha=\Gamma^1_{21}\,\alpha+
\Gamma^2_{21}\,\gamma,
&&\nabla_{\gamma}\gamma=\Gamma^1_{22}\,\alpha+
\Gamma^2_{22}\,\gamma.
\endxalignat
$$
The use of these expansions is correct, since in the first
case of intermediate degeneration the fields $\alpha$ and 
$\gamma$ are non-collinear and they form moving frame
on the plane. Let's denote by $C$ and $D$ the components
of the field $\gamma$:
$$
\xalignat 2
&C=\gamma^1,&&D=\gamma^2.
\tag5.3
\endxalignat
$$
In terms of \thetag{5.3} for the coefficient $\Gamma^1_{22}$
in \thetag{5.2} one can derive
$$
\aligned
\Gamma^1_{22}=&\frac{C\,D\,(C_{\sssize 1.0}-D_{\sssize 0.1})}{M}
+\frac{D^2\,C_{\sssize 0.1}-C^2\,D_{\sssize 1.0}}{M}+\\
\vspace{1ex}
&\quad+\frac{P\,C^3+3\,Q\,C^2\,D+3\,R\,C\,D^2+S\,D^3}{M}.
\endaligned
\tag5.4
$$
The quantity $\Gamma^1_{22}$ given by \thetag{5.4} is a
pseudoscalar field of the weight $4$. According to \cite{24},
it defines the third basic invariant for the first case of
intermediate degeneration
$$
I_3=\frac{\Gamma^1_{22}\,Q^2}{M}.
\tag5.5
$$
By differentiating invariants $I_1$, $I_2$, and $I_3$ along 
pseudovectorial fields $\alpha$ and $\gamma$ we get six new 
invariants $I_5$, $I_6$, $I_7$, $I_8$, $I_9$, and $I_{10}$:
$$
\xalignat 3
&I_4=\frac{\nabla_\alpha  I_1}{N},
&&I_5=\frac{\nabla_\alpha I_2}{N},
&&I_6=\frac{\nabla_\alpha I_3}{N},\\
\vspace{1ex}
&I_7=\frac{(\nabla_\gamma  I_1)^2}{N^3},
&&I_8=\frac{(\nabla_\gamma  I_2)^2}{N^3},
&&I_9=\frac{(\nabla_\gamma  I_3)^2}{N^3}.
\endxalignat
$$
Repeating this procedure more and more, we can form an indefinite
sequence of scalar invariants $I_1$, $I_2$, $I_3$, \dots, adding
$6$ ones in each step.\par
    The number of basic invariants in first case of intermediate
degeneration is less by one than in case of general position,
these are the invariants $I_1$ and $I_2$ from \thetag{5.1} and
the invariant $I_3$ from \thetag{5.5}. Other coefficients in
\thetag{5.4} do not change the number of basic invariants, since
values of them are trivial in most:
$$
\xalignat 2
&\Gamma^2_{11}=\Gamma^1_{21}=0,
&&\Gamma^1_{11}=\Gamma^2_{21}=-\frac{3}{5}\,N.
\endxalignat
$$
For nontrivial coefficients $\Gamma^2_{22}=-\Gamma^1_{12}$ and
$\Gamma^2_{12}$ the following relations were derived:
$$
\gather
I_1\,\Gamma^2_{12}=I_4\,N-\frac{3}{5}\,I_1\,N-2\,I_1^2\,N,
\tag5.6\\
\vspace{2ex}
\aligned
 \bigl(I_1\,\Gamma^2_{22}\bigr)^4&+\bigl(I_7\,N^3\bigr)^2
 +\bigl(16\,I_2\,N^3\,{I_1}^4\bigr)^2=\\
 \vspace{1ex}
 &=32\,I_7\,N^6\,I_2\,{I_1}^4
 +2\,\bigl(I_7\,N^3+16\,I_2\,N^3\,{I_1}^4\bigr)\,
 \bigl(I_1\,\Gamma^2_{22}\bigr)^2.\hskip-2em
\endaligned
\tag5.7
\endgather
$$
Due to the relations \thetag{5.6} and \thetag{5.7} derived
in \cite{24} the coefficients $\Gamma^2_{22}$, $\Gamma^1_{12}$, 
and $\Gamma^2_{12}$ can be expressed through the invariants
from the above sequence and the field $N$.\par
     Here, like in the case of general position, the structure
of invariants in the sequence $I_1$, $I_2$, $I_3$, \dots 
distinguishes three different subcases:
\roster
\item in the infinite sequence of invariants $I_k(x,y)$ one can
      find two functionally independent ones;
\item invariants $I_k(x,y)$ are functionally dependent, but not
      all of them are constants;
\item all invariants in the sequence $I_k(x,y)$ are constants.
\endroster
In the first case the group of point symmetries of the equation
\thetag{1.1} is trivial, in the second case it is one-dimensional,
and in the third case it is two-dimensional. When two-dimensional,
this algebra is Abelian if and only if 
$$
\xalignat 2
&I_1=-\frac{12}{5},
&&I_2=0.
\tag5.8
\endxalignat
$$
This is the result from \cite{24}. Now it is absolutely effective,
since conditions \thetag{5.8} can be tested without transforming
the equation \thetag{1.1} to special coordinates.
\head
6. Second case of intermediate degeneration.
\endhead
     Remember, that for all cases of intermediate degeneration
the parameters $A$ and $B$ do not vanish simultaneously, while
parameter $F$ from \thetag{2.2} is zero identically. The above
first case of intermediate degeneration was distinguished by
the additional condition $M\neq 0$. For the second case of
intermediate degeneration this additional condition is replaced
by the following two relationships:
$$
\xalignat 2
&M=0,&&N\neq 0.
\tag6.1
\endxalignat
$$
Additional condition \thetag{6.1} distinguishes four cases:
second, third, fourth and fifth case of intermediate degeneration.
In each of these four cases one can choose special variables that
realize the conditions \thetag{1.8}, and for which the the
coefficients $P$ and $Q$ in \thetag{1.1} are brought to the form
$$
\xalignat 2
&P=0,&&Q=-\frac{5}{12\,x}.
\tag6.2
\endxalignat
$$
In these variables the third coefficient $R$ in \thetag{1.1} 
also has the special form defined by two arbitrary functions
$r(y)$ and $c(y)$:
$$
R=r(y)+c(y)\,|x|^{-1/4}.
\tag6.3
$$
Second case of intermediate degeneration is distinguished
from third, fourth and fifth cases by additional condition
$$
c(y)\neq 0
\tag6.4
$$
for the function $c(y)$ in \thetag{6.3}.\par
     The conditions \thetag{6.1} are written in arbitrary
variables, they do not need effectivization. Therefore we
are only to make effective the condition \thetag{6.4}.
Let's calculate the pseudoscalar field $\Omega$ in special
coordinates. Taking into account \thetag{6.2}, \thetag{6.3},
and \thetag{1.8}, we use the formula \thetag{4.5}:
$$
\Omega=R_{\sssize 1.0}-2\,Q_{\sssize 0.1}=
-c(y)\,\frac{|x|^{-1/4}}{4\,x}.
\tag6.5
$$
Comparing \thetag{6.4} and \thetag{6.5}, we conclude that
condition \thetag{6.4} can be written as
$$
\Omega\neq 0.
\tag6.6
$$
This condition \thetag{6.6} is an effective form for the 
condition \thetag{6.4}. Being fulfilled in special coordinates,
it remains true in any other coordinates too. This is due to
pseudoscalar rule of transformation for $\Omega$.\par
     Let $F=0$ and let the conditions \thetag{6.1} and \thetag{6.6}
be fulfilled, i\. e\. we are in the second case of intermediate
degeneration. Then at the expense of further specialization of
the choice of variables we can bring the condition \thetag{6.4} to
the form $c(y)=1$. As a result of this the relationship \thetag{6.3}
will have the form
$$
R=r(y)+|x|^{-1/4},
\tag6.7
$$
and for the parameter $S$ we can get the following explicit
expression:
$$
S=\sigma(y)\,|x|^{5/4}-4\,s(y)\,x+\frac{4}{3}\,|x|^2
-12\,\frac{r(y)\,|x|^{7/4}}{x}-4\,\frac{|x|^{3/2}}{x}.
\tag6.8
$$
(for more details see \cite{24}). Here $s(y)$ and $\sigma(y)$
are two arbitrary functions of one variable. The above formulas
\thetag{6.2}, \thetag{6.7}, and \thetag{6.8} define the canonical
form of the equation \thetag{1.1} in the second case of intermediate
degeneration. Algebra of point symmetries for such equation is
described by the following theorem from \cite{24}.
\proclaim{Theorem 6.1} In the second case of intermediate 
degeneration algebra of point symmetries of the equation 
\thetag{1.1} is one-dimensional if and only if parameters
$r(y)$, $s(y)$, and $\sigma(y)$ are identically constant:
$$
\xalignat 3
&\quad r'(y)=0,&&s'(y)=0,&&\sigma'(y)=0.
\hskip-2em
\tag6.9
\endxalignat
$$
If at least one of the conditions \thetag{6.9} fails, then
corresponding algebra of point symmetries is trivial.
\endproclaim
     Conditions \thetag{6.9} determining the structure of the
algebra of point symmetries in the theorem~6.1 are written in
special variables. Therefore they require effectivization.
Remember, that the condition $M\neq 0$ is equivalent to the
non-collinearity of pseudovectorial fields $\alpha$ and $\gamma$.
Conversely, from $M=0$ we derive $\gamma\parallel\alpha$ and remember
that $\alpha\neq 0$. Then the proportionality factor relating these
two field $\gamma=\Lambda\,\alpha$ defines one more pseudoscalar
field of the weight $1$. The field $\Lambda$ can be calculated
by one of the following formulas similar to \thetag{4.3}:
$$
\xalignat 2
&\Lambda=\frac{C}{B},&&\Lambda=-\frac{D}{A}.
\tag6.10
\endxalignat
$$
Here by $C$ and $D$ the components of the field $\gamma$ are denoted
(see \thetag{5.3}). Formulas \thetag{6.10} are effective, they are
applicable in arbitrary coordinates. It would be convenient to write
them in more explicit form:
$$
\align
&\Lambda=-\frac{6\,N\,(A\,S-B_{\sssize 0.1})}{5\,B^2}
-\frac{N_{\sssize 0.1}}{B}+\frac{6\,N\,R}{5\,B}-2\,\Omega,
\tag6.11\\
\vspace{2ex}
&\Lambda=\frac{6\,N\,(B\,P+A_{\sssize 1.0})}{5\,A^2}
-\frac{N_{\sssize 1.0}}{A}-\frac{6\,N\,Q}{5\,A}-2\,\Omega.
\tag6.12
\endalign
$$
These two formulas  \thetag{6.11} and \thetag{6.12} are easily
derived from \thetag{4.30} and \thetag{4.33}. They are mirror
symmetric to each other, while mirror transformation for $\Lambda$
is written as follows: $\Lambda\to-\Lambda$.\par
     Let's calculate the fields $\Omega$ and $\Lambda$ in special
coordinates that were introduced above for second case of special
degeneration. From \thetag{6.2}, \thetag{6.7}, and \thetag{6.8} we 
derive
$$
\xalignat 2
&\Omega=-\frac{|x|^{-1/4}}{4\,x},
&&\Lambda=-\frac{r(y)}{2\,x}.
\tag6.13
\endxalignat
$$
By means of $N$, $\Omega$, and $\Lambda$ let's construct one more field
$$
I_1=\frac{\Lambda^{12}}{\Omega^8\,N^2}.
\tag6.14
$$
The weight of the field \thetag{6.14} appears to be zero: $12-8-2\cdot 
2=0$. Thus, the field $I_1$ is a scalar invariant of the equation
\thetag{1.1}. It is not difficult to calculate this field in special 
coordinates:
$$
I_1=\frac{2304}{25}\,r(y)^{12}.
\tag6.15
$$
By comparing \thetag{6.15} with \thetag{6.9} we conclude that
first of the conditions \thetag{6.9} can be written in the following
invariant form:
$$
I_1=\const.
\tag6.16
$$
The condition \thetag{6.16} can be checked effectively without
transforming the equation \thetag{1.1} to special coordinates.
Effectivization of the rest two conditions in \thetag{6.9}
requires some additional efforts.\par
\head
7. Curvature tensor and additional fields.
\endhead
     Let $F=0$ and $M=0$, while parameters $A$ and $B$ do not
vanish simultaneously. This corresponds to any case of intermediate 
degeneration, except for the first. Theorems~4.1 and 4.2 give us
effective formulas for the parameters $\varphi_1$ and $\varphi_2$.
Then these parameters are used to determine the connection
components \thetag{4.22}. In turn they determine the field of
curvature tensor:
$$
R^k_{qij}=
\frac{\partial\varGamma^k_{jq}}{\partial u^i}
-\frac{\partial\varGamma^k_{iq}}{\partial u^j}+
\sum^2_{s=1}\varGamma^k_{is}\varGamma^s_{jq}-
\sum^2_{s=1}\varGamma^k_{js}\varGamma^s_{iq}.
\tag7.1
$$
Curvature tensor \thetag{7.1} is skew symmetric with respect to the
last pair of indices $i$ and $j$. I two-dimensional geometry such
tensor can be decomposed as $R^k_{qij}=R^k_q\,d_{ij}$. Here $R^k_q$
is pseudotensorial field of the weight $1$. It can be calculated
by the formula
$$
R^k_q=\frac{1}{2}\sum^2_{i=1}\sum^2_{j=1}R^k_{qij}\,d^{ij}.
\tag7.2
$$\par
     For to study the property of the pseudotensorial field
\thetag{7.2} let's calculate its components in special coordinates,
where the conditions \thetag{1.8} hold. In such coordinates the
field $\alpha$ has unitary components: $\alpha^1=1$, $\alpha^2=0$.
The conditions $F=0$ and $M=0$ are written as $P=0$ and 
$Q_{\sssize 1.0}=5/12\,Q^2$, and the parameters $\varphi_1$ and
$\varphi_2$ are given by the formulas \thetag{4.7}. Taking into
account all these circumstances, we can calculate the components
of the field \thetag{7.2} in explicit form. It is remarkable that
due to $M=0$ the matrix $R^k_q$ appears to be upper-triangular
in special coordinates:
$$
R^k_q=\Vmatrix
R^1_1 & R^1_2\\
\vspace{1ex}
0     & R^2_2
\endVmatrix.
\tag7.3
$$
Eigenvalues $\lambda_1=R^1_1$ and $\lambda_2=R^2_2$ of the
matrix \thetag{7.3} are the pseudoscalar fields of the
weight $1$. They can be calculated in explicit form:
$$
\xalignat 2
&\quad\lambda_1=-\frac{3}{5}\,\Lambda,
&&\lambda_2=\frac{3}{5}\,\Omega+\frac{3}{5}\,\Lambda.
\tag7.4
\endxalignat
$$
Let subtract the identity matrix multiplied by the second 
eigenvalue $\lambda_2=R^2_2$ from the matrix \thetag{7.3}. 
As a result we have the matrix
$$
P^k_q=R^k_q-\lambda_2\,\delta^k_q=\Vmatrix
(\lambda_1-\lambda_2) & R^1_2\\
\vspace{4ex}
0     &    0
\endVmatrix,
\tag7.5
$$
which defines another pseudotensorial field of the weight $1$.
Let $\bold X$ be some arbitrary vector-field with components
$X^1$ and $X^2$. Let's contract it with \thetag{7.5}. Then
we obtain pseudovectorial field $P\bold X$ of the weight $1$,
whose second component being identically zero. This means that
$P\bold X$ is collinear to the field $\alpha$ having unitary
components $\alpha^1=1$ and $\alpha^2=0$ in special coordinates:
$$
P\bold X=\left((\lambda_1-\lambda_2)\,X^1+R^1_2\,X^2\right)
\,\alpha.
\tag7.6
$$
The proportionality factor binding $P\bold X$ and $\alpha$ in
\thetag{7.6} depends linearly on the components of the vector
$\bold X$. Therefore it defines pseudocovectorial field of the
weight $-1$ with the following components:
$$
\xalignat 2
&\omega_1=\lambda_1-\lambda_2,
&&\omega_2=R^1_2.
\tag7.7
\endxalignat
$$
Components $\omega_1$ and $\omega_2$ from \thetag{7.7} can be
calculated explicitly:
$$
\xalignat 2
&\quad\omega_1=-\frac{3}{5}\,\Omega-\frac{6}{5}\,\Lambda,
&&\omega_2=S_{\sssize 1.0}-\frac{6}{5}\,R_{\sssize 0.1}
+\frac{12}{5}\,S\,Q-\frac{54}{25}\,R^2.
\hskip-3em
\tag7.8
\endxalignat
$$\par
     Formulas \thetag{7.4} for the eigenvalues of the field
\thetag{7.3} need no effectivization. They do not change 
in arbitrary coordinates if we take into account the formulas
\thetag{4.17}, \thetag{4.21}, \thetag{6.11}, and \thetag{6.12} 
for the fields $\Omega$ and $\Lambda$. But the formulas
\thetag{7.8} should be recalculated for the case of arbitrary
coordinates. In order to do it note that the formulas \thetag{7.1} 
and \thetag{7.2} hold for arbitrary coordinates. Then formula
\thetag{7.5} in non-special coordinates is written as
$$
P^k_q=R^k_q-\lambda_2\,\delta^k_q=
\Vmatrix
R^1_1-\lambda_2 & R^1_2\\
\vspace{4ex}
R^2_1  &    R^2_2-\lambda_2
\endVmatrix.
$$
This is due to the fact that matrix $R^k_q$ in arbitrary coordinates
isn't upper-triangular. However, the collinearity of the fields
$P\bold X$ and $\alpha$ does not depend on the choice of coordinates.
Hence the formula \thetag{7.6} in arbitrary coordinates is written as
$$
P\bold X=
\Vmatrix
R^1_1-\lambda_2 & R^1_2\\
\vspace{4ex}
R^2_1  &    R^2_2-\lambda_2
\endVmatrix\cdot
\Vmatrix X^1\\ \vspace{4ex} X^2\endVmatrix=
(\omega_1\,X^1+\omega_2\,X^2)\cdot
\Vmatrix B\vphantom{X^1}\\ \vspace{4ex} -A\endVmatrix.
\hskip-3em
\tag7.9
$$
From \thetag{7.9} we easily extract the required effective
formulas for $\omega_1$ and $\omega_2$:
$$
\xalignat 2
&\omega_1=\frac{R^1_1-\lambda_2}{B},
&\omega_2=\frac{R^1_2}{B}.
\tag7.10
\endxalignat
$$
Formulas \thetag{7.10} hold for $B\neq 0$. For the case $A\neq 0$
we can write mirror symmetric formulas. They are the following
ones:
$$
\xalignat 2
&\omega_1=-\frac{R^2_1}{A},
&&\omega_2=\frac{\lambda_2-R^2_2}{A}.
\tag7.11
\endxalignat
$$
Upon explicit calculation of the components of matrix $R^k_q$
and upon substituting them into \thetag{7.10} for 
the field $\omega$ in arbitrary coordinates we get
$$
\gather
\gathered
\omega_1=-\frac{6\,\Lambda+3\,\Omega}{5\,B}
+\frac{5\,A\,S_{\sssize 1.0}-6\,A\,R_{\sssize 0.1}+
12\,Q\,A\,S}{5\,B^2}-\frac{54}{25}\,\frac{A\,R^2}{B^2}+\\
\vspace{1ex}
+\frac{2\,A\,A_{\sssize 0.1}\,S+A\,B_{\sssize 1.0}\,S+
A^2\,S_{\sssize 0.1}-A\,B_{\sssize 0.2}}{5\,B^3}
-\frac{12\,A^2\,S\,R}{25\,B^3}+\\
\vspace{1ex}
+\frac{3\,A\,R\,B_{\sssize 0.1}}{25\,B^3}
+\frac{6\,A\,B_{\sssize 0.1}^2
+6\,A^3\,S^2-12\,A^2\,B_{\sssize 0.1}\,S}{25\,B^4},
\endgathered
\tag7.12\\
\vspace{2ex}
\gathered
\omega_2=\frac{12\,S\,Q}{5\,B}-\frac{54}{25}\,\frac{R^2}{B}+
\frac{S_{\sssize 1.0}}{B}-\frac{6\,R_{\sssize 0.1}}{5\,B}
+\frac{S\,B_{\sssize 1.0}+A\,S_{\sssize 0.1}
-B_{\sssize 0.2}}{5\,B^2}+\\
\vspace{1ex}
\quad+\frac{2\,A_{\sssize 0.1}\,S}{5\,B^2}
-\frac{3\,R\,B_{\sssize 0.1}+12\,S\,A\,R}{25\,B^2}
+\frac{6\,A^2\,S^2-12\,B_{\sssize 0.1}\,A\,S
+6\,B_{\sssize 0.1}^2}{25\,B^3}.
\hskip-3em
\endgathered
\hskip-2em
\tag7.13
\endgather
$$
Formulas \thetag{7.12} and \thetag{7.13} hold for $B\neq 0$.
For the case $A\neq 0$ we have formulas that are derived
from \thetag{7.11}:
$$
\gather
\gathered
\omega_1=\frac{12\,P\,R}{5\,A}-\frac{54}{25}\,\frac{Q^2}{A}-
\frac{P_{\sssize 0.1}}{A}+\frac{6\,Q_{\sssize 1.0}}{5\,A}
-\frac{P\,A_{\sssize 0.1}+B\,P_{\sssize 1.0}
+A_{\sssize 2.0}}{5\,A^2}-\\
\vspace{1ex}
\quad-\frac{2\,P\,B_{\sssize 1.0}}{5\,A^2}
+\frac{3\,Q\,A_{\sssize 1.0}-12\,P\,B\,Q}{25\,A^2}
+\frac{6\,B^2\,P^2+12\,B\,P\,A_{\sssize 1.0}
+6\,A_{\sssize 1.0}^2}{25\,A^3},
\hskip-3em
\endgathered
\hskip-2em
\tag7.14\\
\vspace{2ex}
\gathered
\omega_2=\frac{6\,\Lambda+3\,\Omega}{5\,A}
+\frac{6\,B\,Q_{\sssize 1.0}+12\,R\,B\,P
-5\,B\,P_{\sssize 0.1}}{5\,A^2}
-\frac{54}{25}\,\frac{B\,Q^2}{A^2}-\\
\vspace{1ex}
-\frac{2\,B\,B_{\sssize 1.0}\,P+B\,P\,A_{\sssize 0.1}+
B^2\,P_{\sssize 1.0}+B\,A_{\sssize 2.0}}{5\,A^3}
-\frac{12\,B^2\,P\,Q}{25\,A^3}+\\
\vspace{1ex}
+\frac{3\,B\,Q\,A_{\sssize 1.0}}{25\,A^3}
+\frac{6\,B\,A_{\sssize 1.0}^2
+6\,B^3\,P^2+12\,B^2\,P\,A_{\sssize 1.0}}{25\,A^4}.
\endgathered
\tag7.15
\endgather
$$
Formulas \thetag{7.14} and \thetag{7.15} are mirror symmetric
with respect to the formulas \thetag{7.12} and \thetag{7.13}. 
Comparing these two pairs of formulas, we derive the following 
mirror transformations for $\omega_1$ and $\omega_2$:
$$
\xalignat 2
&\omega_1\to-\omega_2,
&&\omega_2\to-\omega_1.
\tag7.16
\endxalignat
$$
Mirror transformations \thetag{7.16} are to be considered as
the expansion for the transformations \thetag{4.18}, \thetag{4.19}, 
\thetag{4.20}, and \thetag{4.34}.\par
     Let $F=0$ and $M=0$ as before and suppose that $A$ and $B$
do not vanish simultaneously. Under these assumptions we have two
pseudoscalar fields $\Omega$ and $\Lambda$ of the weight $1$.
In special coordinates they are defined by
$$
\xalignat 2
&\Omega=R_{\sssize 1.0}-2\,Q_{\sssize 0.1},
&&\Lambda=-2\,R_{\sssize 1.0}+3\,Q_{\sssize 0.1}+\frac{6}{5}\,Q\,R.
\tag7.17
\endxalignat
$$
The effectivization for the formulas \thetag{7.17} has been already 
done in form of relationships \thetag{4.17}, \thetag{4.21}, 
\thetag{6.11}, and \thetag{6.12}. Let's calculate the covariant
differentials $\nabla\Omega$ and $\nabla\Lambda$ for $\Omega$ and
$\Lambda$. These are pseudocovectorial fields of the weight $1$.
From \thetag{7.17} and \thetag{4.23} for the components of
$\nabla\Omega$ and $\nabla\Lambda$ in special coordinates we
derive
$$
\xalignat 2
&\nabla_1\Omega=\frac{9}{5}\,Q\,\Omega,
&&\nabla_2\Omega=\Omega_{\sssize 0.1}-\frac{3}{5}\,R\,\Omega,
\tag7.18\\
\allowdisplaybreak
&\nabla_1\Lambda=\frac{6}{5}\,Q\,\Lambda,
&&\nabla_2\Lambda=\Lambda_{\sssize 0.1}-\frac{3}{5}\,R\,\Lambda.
\tag7.19
\endxalignat
$$
Let's compare \thetag{7.18} and \thetag{7.19} with the formulas
\thetag{7.8} for the components of the field $\omega$, which
also has the weight $1$. Then let's construct the field
$$
w=N\,\omega+\nabla\Lambda+\frac{1}{3}\,\nabla\Omega.
\tag7.20
$$
This field \thetag{7.20} has the weight $1$. It's remarkable
that its first component in special coordinates is zero:
$w_1=0$. Hence $w$ is collinear to pseudocovectorial field
$\alpha$ with components $\alpha_1=0$ and $\alpha_2=1$. Let's
denote by $K$ the proportionality factor in $w=K\,\alpha$.
Then $K$ is a scalar field (field of the weight $0$). In special
coordinates it is calculated as follows:
$$
\aligned
K&=\Lambda_{\sssize 0.1}-\frac{3}{5}\,R\,\Lambda
+\frac{1}{3}\,\Omega_{\sssize 0.1}-\frac{1}{5}\,R\,\Omega+\\
\vspace{1ex}
&+Q\,S_{\sssize 1.0}-\frac{6}{5}\,Q\,R_{\sssize 0.1}
+\frac{12}{5}\,S\,Q^2-\frac{54}{25}\,Q\,R^2.
\endaligned
\tag7.21
$$
It's no problem to make effective \thetag{7.21}, since the field
\thetag{7.20} can be evaluated in arbitrary coordinates. When
$B\neq 0$, we have
$$
K=\frac{\Lambda_{\sssize 0.1}+\Lambda\varphi_2}{B}+
\frac{\Omega_{\sssize 0.1}+\Omega\varphi_2}{3\,B}+
\frac{N\,\omega_2}{B}.
\tag7.22
$$
Here $\Lambda$ and $\Omega$ are calculated by the formulas
\thetag{6.11} and \thetag{4.17}, $N$ is defined by the first
relationship \thetag{4.3}, parameter $\varphi_2$ is given by
the formula \thetag{4.8}, and $\omega_2$ is defined by the
formula \thetag{7.13}. For $A\neq 0$ we have the formula mirror 
symmetric to \thetag{7.22}:
$$
K=\frac{\Lambda_{\sssize 1.0}+\Lambda\varphi_1}{A}+
\frac{\Omega_{\sssize 1.0}+\Omega\varphi_1}{3\,A}+
\frac{N\,\omega_1}{A}.
\tag7.23
$$
In \thetag{7.23} fields $\Lambda$ and $\Omega$ are calculated
by \thetag{6.12} and \thetag{4.21}, field $N$ is defined by the
second relationship \thetag{4.3}, parameter $\varphi_1$ is given
by \thetag{4.14}, and $\omega_1$ is calculated by the formula
\thetag{7.14}.
\head
8. Algebra of symmetries in the second case of intermediate
degeneration.
\endhead
     In special coordinates the structure of the algebra of point
symmetries of the equation \thetag{1.1} for this case is described
by the conditions \thetag{6.9} in theorem~6.1. One of them had been
made effective in form of the condition \thetag{6.16}. In order to 
make effective two other conditions \thetag{6.9} we shall construct
some additional scalar invariants. Let's consider pseudocovectorial
field $\varepsilon$ defined by the formula analogous to \thetag{7.20}:
$$
\varepsilon=N\,\omega+\nabla\Lambda.
\tag8.1
$$
Field \thetag{8.1} has the weight $1$. After raising indices by means
of skew-symmetric matrix from \thetag{2.5} we get the pseudovectorial
field $\varepsilon$ of the weight $2$. It's easy to calculate the
components of this field in special coordinates:
$$
\aligned
\varepsilon^1&=Q\,S_{\sssize 1.0}-\frac{6}{5}\,Q\,R_{\sssize 0.1}
+\frac{12}{5}\,S\,Q^2-\frac{54}{25}\,Q\,R^2+\Lambda_{\sssize 0.1}
-\frac{3}{5}\,R\,\Lambda,\\
\varepsilon^2&=\frac{3}{5}\,Q\,\Omega.
\endaligned
\tag8.2
$$
In the second case of intermediate degeneration we have $N\neq 0$ 
and $\Omega\neq 0$ (see \thetag{6.1} and \thetag{6.6}).
In special coordinates this yields $Q\,\Omega\neq 0$, i\. e\.
second component of the field \thetag{8.2} is nonzero.
Therefore pseudovectorial fields $\alpha$ and $\varepsilon$ 
of the same weight $2$ are non-collinear and we are able
to write the expansions
$$
\xalignat 2
&\nabla_{\alpha}\alpha=\Gamma^1_{11}\,\alpha+
\Gamma^2_{11}\,\varepsilon,
&&\nabla_{\alpha}\varepsilon=\Gamma^1_{12}\,\alpha+
\Gamma^2_{12}\,\varepsilon,\\
\vspace{-1.7ex}
&&&\tag8.3\\
\vspace{-1.7ex}
&\nabla_{\varepsilon}\alpha=\Gamma^1_{21}\,\alpha+
\Gamma^2_{21}\,\varepsilon,
&&\nabla_{\varepsilon}\varepsilon=\Gamma^1_{22}\,\alpha+
\Gamma^2_{22}\,\varepsilon.
\endxalignat
$$
Covariant derivatives in \thetag{8.3} are determined by the
connection \thetag{4.22}. All coefficients $\Gamma^k_{ij}$ in
these expansions are the pseudoscalar fields of the weight $2$.
Most of them are trivial or equal to zero:
$$
\xalignat 2
&\Gamma^2_{11}=\Gamma^1_{21}=0,
&&\Gamma^1_{11}=\Gamma^2_{21}=-\frac{3}{5}\,N.
\endxalignat
$$
Others are less trivial, but, nevertheless, they can be expressed
through the pseudoscalar fields $K$, $N$, $\Omega$, and $\Lambda$:
$$
\aligned
\Gamma^1_{12}&=-\frac{6}{5}\,K\,N+N-\frac{6}{5}\,\Lambda^2-
3\,\Omega\,\Lambda,\\
\vspace{1ex}
\Gamma^2_{22}&=\frac{9}{5}\,K\,N-\frac{6}{5}\,\Omega^2
-\frac{3}{5}\,\Omega\,\Lambda.
\endaligned
$$
The only new field in \thetag{8.3} is the field $\Gamma^1_{22}$.
Explicit formula for $\Gamma^1_{22}$ in special coordinates 
contains $56$ summands. We shall not write it here. Instead of this,
we shall describe the effective algorithm to calculate $\Gamma^1_{22}$ 
in arbitrary coordinates.\par
     Let's start with the relationships \thetag{8.2} for the components
of pseudovectorial field $\varepsilon$. Their effectivization is based
on the formula \thetag{8.1} for this field:
$$
\xalignat 2
&\quad\varepsilon^1=N\,\omega_2+\Omega_{\sssize 0.1}
+\varphi_2\,\Omega,
&&\varepsilon^2=-N\,\omega_1-\Omega_{\sssize 1.0}
+\varphi_1\,\Omega.
\hskip-3em
\tag8.4
\endxalignat
$$
Here $\omega_1$ and $\omega_2$ are calculated by the formulas
\thetag{7.12}, \thetag{7.13}, \thetag{7.14}, and \thetag{7.15},
parameters $\varphi_1$ and $\varphi_2$ are given by \thetag{4.8}
and \thetag{4.14}. Denote by $C$ and $D$ the components of the
field \thetag{8.4}: $\varepsilon^1=C$ and $\varepsilon^2=D$. 
Then
$$
\aligned
\Gamma^1_{22}=&\frac{5\,D\,C\,(C_{\sssize 1.0}
-D_{\sssize 0.1})}{3\,N\,\Omega}+\frac{5\,D^2\,C_{\sssize 0.1}
-5\,C^2\,D_{\sssize 1.0}}{3\,N\,\Omega}+\\
\vspace{1ex}
&\qquad+\frac{5\,P\,C^3+15\,Q\,C^2\,D+15\,R\,C\,D^2+5\,S\,D^3}
{3\,N\,\Omega}.
\endaligned
\hskip-3em
\tag8.5
$$
Let's calculate the field $K$ in special coordinates, where
the coefficients of the equation \thetag{1.1} is defined by
the formulas \thetag{6.2}, \thetag{6.7}, and \thetag{6.8}:
$$
K=-\frac{5}{48}\,\sigma(y)\,|x|^{-3/4}-\frac{5}{9}
+\frac{9}{10}\,\frac{r(y)\,|x|^{-1/4}}{x}
+\frac{7}{60}\,\frac{|x|^{-1/2}}{x}+\frac{6\,r(y)^2}{5\,x}.
$$
Relying on this formula and on the relationships \thetag{6.13},
we shall construct the following two pseudoscalar fields:
$$
\xalignat 2
&L=K\,N+\frac{5}{9}\,N+3\,\Lambda\,\Omega+\frac{7}{9}\,\Omega^2
+2\,\Lambda^2,
&&I_2=\frac{L^4}{N^2\,\Omega^4}.
\tag8.6
\endxalignat
$$
Field $L$ has the weight $2$, while the weight of the field $I_2$
is zero, i\. e\. $I_2$ is a scalar invariant of the equation
\thetag{1.1}. Fields \thetag{8.6} can be evaluated explicitly:
$$
\xalignat 2
&L=\frac{25}{576}\,\frac{\sigma(y)\,|x|^{-3/4}}{x},
&&I_2=\frac{15625}{2985984}\,\sigma(y)^4.
\tag8.7
\endxalignat
$$
Second formula in \thetag{8.7} is similar to \thetag{6.15}.
Due to this formula we can make effective the second condition
\thetag{6.9} from theorem~6.1:
$$
I_2=\const.
\tag8.8
$$\par
     Now we have only to make effective the rest third condition 
in \thetag{6.9}. First, we evaluate the covariant derivatives of
$\Lambda$ and $L$ along the pseudovectorial field $\varepsilon$.
The field $\nabla_\varepsilon\Lambda$ has the weight $3$, the weight
of $\nabla_\varepsilon L$ is $4$. Then we combine these two fields
with the field \thetag{8.5}, which has the weight $2$:
$$
\align
E&=\Gamma^1_{22}-\frac{\nabla_\varepsilon L}{N}
+\frac{4\,\Lambda\,\nabla_\varepsilon\Lambda}{N}
+\frac{17\,\Omega\,\nabla_\varepsilon\Lambda}{6\,N}
+\frac{12\,L^2}{5\,N}
-\frac{53\,L\,\Lambda\,\Omega}{5\,N}-\\
\vspace{1ex}
&-\frac{48\,L\,\Lambda^2}{5\,N}
-\frac{62\,L\,\Omega^2}{15\,N}
-\frac{8\,L}{3}
+\frac{48\,\Lambda^4}{5\,N}
+\frac{106\,\Lambda^3\,\Omega}{5\,N}
+\frac{16\,\Lambda^2}{3}+
\tag8.9\\
\vspace{1ex}
&+\frac{1163\,\Lambda^2\,\Omega^2}{60\,N}
+\frac{137\,\Lambda\,\Omega^3}{18\,N}
+\frac{50\,\Lambda\,\Omega}{9}
+\frac{203\,\Omega^2}{108}
+\frac{77\,\Omega^4}{135\,N}
+\frac{20\,N}{27}.
\endalign
$$
Field $E$ in \thetag{8.9} has the weight $2$. It is intentionally
constructed so that its value in special coordinates is proportional
to the function $s(y)$ from \thetag{6.8}:
$$
E=-\frac{s(y)\,|x|^{-1/2}}{64\,x^3}.
\tag8.10
$$
Due to this relationship \thetag{8.10} we can build the third
scalar invariant for the second case of intermediate degeneration:
$$
I_3=\frac{E^6\,N^4}{\Omega^{20}}.
\tag8.11
$$
The value of the invariant \thetag{8.11}  in special coordinates
doesn't depend on $x$:
$$
I_3=\frac{625}{1296}\,s(y)^6.
$$
This relationship is used to write the third condition \thetag{6.9} 
in invariant form:
$$
I_3=\const.
\tag8.12
$$
Now the theorem~6.1 can be reformulated as follows.
\proclaim{Theorem 8.1} In the second case of intermediate 
degeneration algebra of point symmetries of the equation 
\thetag{1.1} is one-dimensional if and only if the conditions
\thetag{6.16}, \thetag{8.8}, and \thetag{8.12} hold, i\. e\. 
if all invariants $I_1$, $I_2$, and $I_3$ are identically
constant. If at least one of these conditions fails, then 
corresponding algebra of point symmetries is trivial.
\endproclaim
\head
9. Third case of intermediate degeneration.
\endhead
     In third case of intermediate degeneration we have
$F=0$. Parameters $A$ and $B$ do not vanish simultaneously.
The conditions \thetag{6.1} are fulfilled too. But the
condition \thetag{6.4} from the second case of intermediate
degeneration is replaced by the following pair of
relationship written in special coordinates:
$$
\xalignat 2
&c(y)=0,&&r(y)\neq 0.
\tag9.1
\endxalignat
$$
First of these conditions is exactly opposite to \thetag{6.4}.
In effective form it is written as $\Omega=0$ (compare with
the condition \thetag{6.6} above). By means of direct 
calculations for the field $\Lambda$ we derive 
$$
\Lambda=-\frac{r(y)}{2\,x}.
\tag9.2
$$
Due to \thetag{9.2} the conditions \thetag{9.1} distinguishing 
third case of intermediate degeneration are written in the
following invariant form:
$$
\xalignat 2
&\Omega=0,&&\Lambda\neq 0.
\tag9.3
\endxalignat
$$
According to the results of \cite{24}, in this case one can
find the special coordinates such that the coefficients
$P$, $Q$, and $R$ of the equation \thetag{1.1} are brought
to the form
$$
\xalignat 3
&\quad P=0,&&Q=-\frac{5}{12\,x},&&R=1.
\hskip-2em
\tag9.4
\endxalignat
$$
The fourth coefficient $S$ in the equation \thetag{1.1} is
brought to the form
$$
S=\sigma(y)\,|x|^{5/4}-4\,s(y)\,x
+\frac{4}{3}\,|x|^2.
\tag9.5
$$
From \thetag{9.4} and \thetag{9.5} one can derive the following
theorem (see \cite{24}).
\proclaim{Theorem 9.1} In the third case of intermediate 
degeneration algebra of point symmetries of the equation 
\thetag{1.1} is one-dimensional if and only if parameters
$s(y)$ and $\sigma(y)$ in \thetag{9.5} are identically 
constant:
$$
\xalignat 3
&s'(y)=0,&&\sigma'(y)=0.
\hskip-2em
\tag9.6
\endxalignat
$$
If at least one of the conditions \thetag{9.6} fails, then
corresponding algebra of point symmetries is trivial.
\endproclaim    
    In order to make effective \thetag{9.6} let's consider again
the field $\omega$ from \thetag{7.8}. This field has the weight $-1$. 
Upon raising indices by means of the matrix \thetag{2.5} we get the
vector field $\omega$. In special coordinates its components are 
the following:
$$
\xalignat 2
&\quad\omega^1=S_{\sssize 1.0}-\frac{6}{5}\,R_{\sssize 0.1}
+\frac{12}{5}\,S\,Q-\frac{54}{25}\,R^2,
&&\omega^2=\frac{6}{5}\,\Lambda.
\hskip-3em
\tag9.7
\endxalignat
$$
Components of the vector-fields \thetag{9.7} in arbitrary
coordinates can be effectively calculated by the formulas
\thetag{7.12}, \thetag{7.13}, \thetag{7.14}, and \thetag{7.15}.
One should only take into account that $\omega^1=\omega_2$ and 
$\omega^2=-\omega_1$.\par
     From \thetag{9.7} we see that the condition of non-collinearity
of the fields $\omega$ and $\alpha$  coincides with $\Lambda\neq 0$.
In third case of intermediate degeneration this condition holds
(see \thetag{9.3}). We shall use it to write the following
expansions analogous to \thetag{8.3}:
$$
\xalignat 2
&\nabla_{\alpha}\alpha=\Gamma^1_{11}\,\alpha+
\Gamma^2_{11}\,\omega,
&&\nabla_{\alpha}\omega=\Gamma^1_{12}\,\alpha+
\Gamma^2_{12}\,\omega,\\
\vspace{-1.7ex}
&&&\tag9.8\\
\vspace{-1.7ex}
&\nabla_{\omega}\alpha=\Gamma^1_{21}\,\alpha+
\Gamma^2_{21}\,\omega,
&&\nabla_{\omega}\omega=\Gamma^1_{22}\,\alpha+
\Gamma^2_{22}\,\omega,
\endxalignat
$$
The expansions \thetag{9.8} define the series of pseudoscalar
fields: $\Gamma^2_{11}$ is the field of the weight $4$, fields
$\Gamma^1_{11}$, $\Gamma^2_{12}$, and $\Gamma^2_{21}$ have the
weight $2$, next three fields $\Gamma^1_{12}$, $\Gamma^1_{21}$, 
and $\Gamma^2_{22}$ have the weight $0$, and the last field
$\Gamma^1_{22}$ is of the weight $-2$. Most of these fields
can be reduced to the various combinations of the fields that
were already defined:
$$
\xalignat 2
&\Gamma^2_{11}=\Gamma^1_{21}=0,
&&\Gamma^1_{11}=\Gamma^2_{21}=-\frac{3}{5}\,N.
\endxalignat
$$
The same is true for the next three fields too:
$$
\xalignat 3
&\Gamma^1_{12}=1-\frac{3}{5}\,K,
&&\Gamma^2_{12}=\frac{9}{5}\,N,
&&\Gamma^2_{22}=\frac{6}{5}\,K.
\endxalignat
$$
The only exception is the field $\Gamma^1_{22}$. In special 
coordinates this pseudoscalar field of the weight $-2$ can be 
calculated explicitly:
$$
\aligned
\Gamma^1_{22}&=S_{\sssize 1.0}-\frac{6}{5}\,R_{\sssize 0.1}
-\frac{54}{25}\,R^2+\frac {12}{5}\,S\,Q
+\frac{6}{5}\,\Lambda\,S_{\sssize 1.1}-\\
\vspace{1ex}
&-\frac{36}{25}\,\Lambda\,R_{\sssize 0.2}
-\frac{36}{25}\,\Lambda^2\,S
-\frac{3402}{625}\,R^3\,\Lambda
-\frac{1026}{125}\,R_{\sssize 0.1}\,R\,\Lambda+\\
\vspace{1ex}
&+\frac{72}{25}\,\Lambda\,S_{\sssize 0.1}\,Q
+\frac{63}{25}\,S_{\sssize 1.0}\,R\,\Lambda
-\frac{9}{5}\,S_{\sssize 1.0}\,\Lambda_{\sssize 0.1}
+\frac{54}{25}\,R_{\sssize 0.1}\,\Lambda_{\sssize 0.1}+\\
\vspace{1ex}
&\qquad+\frac{486}{125}\,R^2\,\Lambda_{\sssize 0.1}
+\frac{1188}{125}\,S\,Q\,R\,\Lambda
-\frac{108}{25}\,S\,Q\,\Lambda_{\sssize 0.1}.
\endaligned
\tag9.9
$$
Let's introduce the following notations:
$$
\xalignat 2
&C=\omega^1=\omega_2,
&&D=\omega^2=-\omega_1.
\tag9.10
\endxalignat
$$
Then for $\Gamma^1_{22}$ we can derive the formula similar
to \thetag{5.4} and \thetag{8.5}:
$$
\aligned
\Gamma^1_{22}&=\frac{5\,C\,D\,(C_{\sssize 1.0}
-D_{\sssize 0.1})}{6\,\Lambda}+\frac{5\,D^2\,C_{\sssize 0.1}
-5\,C^2\,D_{\sssize 1.0}}{6\,\Lambda}+\\
\vspace{1ex}
&+\frac{5\,P\,C^3+15\,Q\,C^2\,D+15\,R\,C\,D^2+5\,S\,D^3}
{6\,\Lambda}.
\endaligned
\tag9.11
$$
Formulas \thetag{9.10} and \thetag{9.11} together with
\thetag{7.12}, \thetag{7.13}, \thetag{7.14}, and \thetag{7.15}
give the effective way to calculate the field \thetag{9.9} in
arbitrary coordinates.\par
     Now let's proceed with constructing the invariants for 
the equation \thetag{1.1} in the third case of intermediate 
degeneration. For $L$ and $I_1$ in this case we take
$$
\xalignat 2
&L=K+\frac{5}{9}+\frac{2\,\Lambda^2}{N},
&&I_1=\frac{L^8\,N^6}{\Lambda^{12}}.
\tag9.12
\endxalignat
$$
Both fields $L$ and $I_1$ in \thetag{9.12} have the weigh $0$,
i\. e\. they are scalar invariants. For their values from
\thetag{9.4} and \thetag{9.5} we derive
$$
\xalignat 2
&\quad L=-\frac{5\,|x|^{-3/4}}{48}\,\sigma(y),
&&I_1=\frac{6103515625}{20542695432781824}\,\sigma(y)^8.
\hskip-4em
\tag9.13
\endxalignat
$$
Pseudoscalar field $E$ having the weight $-2$ is composed 
of the fields $\omega$, $L$, $\Omega$, and $N$:
$$
E=\Gamma^1_{22}-\frac{\nabla_\omega L}{N}
 +\frac{9\,L^2}{5\,N}-\frac{2\,L}{N}
 -\frac{12\,L\,\Lambda^2}{5\,N^2}+\frac{7\,\Lambda^2}{3\,N^2}
 +\frac{5}{9\,N}+\frac{63\,\Lambda^4}{20\,N^3}.
\hskip-3em
\tag9.14
$$
We shall use this field \thetag{9.14} in order to construct the
second scalar invariant $I_2$:
$$
I_2=\frac{E\,N^3}{\Lambda^4}.
\tag9.15
$$
Let's evaluate the fields $E$ and $I_2$ using \thetag{9.4},
\thetag{9.5}, \thetag{9.14}, and \thetag{9.15}:
$$
\xalignat 2
&E=-\frac{36}{25\,x}\,s(y),
&&I_2=\frac{5}{3}\,s(y).
\tag9.16
\endxalignat
$$
Formulas \thetag{9.13} and \thetag{9.16} can be used to
reformulate theorem~9.1 in effective form.
\proclaim{Theorem 9.2} In the third case of intermediate 
degeneration algebra of point symmetries of the equation 
\thetag{1.1} is one-dimensional if and only if both
invariants $I_1$ and $I_2$ in \thetag{9.12} and \thetag{9.15}
are identically constant. If at least one of these two conditions 
fails, then corresponding algebra of point symmetries is trivial.
\endproclaim    
\head
10. Fourth case of intermediate degeneration.
\endhead
     In the fourth case of intermediate degeneration we have
$F=0$, and parameters $A$ and $B$ do not vanish simultaneously.
Moreover, we have previous conditions $N\neq 0$ and $M=0$
from \thetag{6.1}, and also new ones
$$
\xalignat 2
&\Omega=0,&&\Lambda=0.
\tag10.1
\endxalignat
$$
All above conditions, including \thetag{10.1}, don't require
special effectivization. They can be tested in arbitrary
coordinates.\par
     According to the results of \cite{24}, when the above
conditions hold, one can choose special coordinates in which
parameters $P$, $Q$, and $R$ in \thetag{1.1} are brought
to the form
$$
\xalignat 3
&\quad P=0,&&Q=-\frac{5}{12\,x},&&R=0.
\hskip-2em
\tag10.2
\endxalignat
$$
The coefficient $S$ in these coordinates also has the special form:
$$
S=\sigma(y)\,|x|^{5/4}
+\frac{4}{3}\,|x|^2.
\tag10.3
$$
Fourth case of intermediate degeneration is distinguished by
the additional condition for the function $\sigma(y)$ in 
\thetag{10.3}:
$$
\sigma(y)\neq 0.
\tag10.4
$$
This condition \thetag{10.4} should be made effective. Fortunately
in this case it doesn't require special efforts. One should only
evaluate the field $K$ in special coordinates, using the relationships
\thetag{10.2} and \thetag{10.3}:
$$
K=-\frac{5\,|x|^{-3/4}}{48}\,\sigma(y)-\frac{5}{9}.
\tag10.5
$$
Due to \thetag{10.5} the condition \thetag{10.4} can be written as
$$
K+\frac{5}{9}\neq 0.
\tag10.6
$$
Note that $K$ is a scalar field, its weight is $0$. Therefore,
being fulfilled in special coordinates, condition \thetag{10.6}
remains true for any other coordinates. Hence it is effective 
form for the condition \thetag{10.4}.\par
\proclaim{Theorem 10.1} In the fourth case of intermediate 
degeneration algebra of point symmetries of the equation 
\thetag{1.1} is one-dimensional if and only if the function
$\sigma(y)$ from \thetag{10.3} is the solution of the following 
differential equation:
$$
\sigma'''-5\,\frac{\sigma''\,\sigma'}{\sigma}
+\frac{40}{9}\,\frac{{\sigma'}^3}{\sigma^2}=0.
\tag10.7
$$
Otherwise if $sigma(y)$ doesn't satisfy the equation \thetag{10.7}, 
then the algebra of point symmetries is trivial.
\endproclaim
    This theorem from \cite{24} gives complete description of
the algebra of point symmetries for the fourth case of intermediate
degeneration. But the condition \thetag{10.7} in its statement
is written for the special coordinates. Now should make it
effective. Let's consider pseudovectorial field $\omega$ from
\thetag{7.8}. This field was constructed by means of curvature tensor
\thetag{7.1} for the case when $F=0$ and $M=0$. Under the conditions 
\thetag{10.1} its first component $\omega_1$ equals to zero in
special coordinates. This means that $\omega\parallel\alpha$, and
we can define new scalar field $\Theta$ by the relationship $\omega=
\Theta\,\alpha$. It is easy to find that $\Theta$ has the 
weight $-2$. Here is the expression for $\Theta$ in special coordinates:
$$
\Theta=S_{\sssize 1.0}-\frac{6}{5}\,R_{\sssize 0.1}
+\frac{12}{5}\,S\,Q-\frac{54}{25}\,R^2.
\tag10.8
$$
Fields $K$ and $\Theta$ are bound with each other by very simple
relationship: $K=N\,\Theta$. For the fourth case of intermediate
degeneration $N\neq 0$, therefore $\Theta$ can be evaluated through
$K$.  However,  when  $N=0$  (in  sixth  and  seventh  cases   of 
intermediate
degeneration), both fields $K$ and $N$ vanish simultaneously.
For this reason it's better to determine $\Theta$ from the
relationship $\omega=\Theta\,\alpha$. Due to this relationship
we can make effective the formula \thetag{10.8} for $\Theta$:
$$
\xalignat 2
&\Theta=\frac{\omega_2}{B},
&&\Theta=\frac{\omega_1}{A}.
\tag10.9
\endxalignat
$$
First of the formulas \thetag{10.9} is used when $B\neq 0$, the
value of $\omega_2$ being calculated by \thetag{7.13}. Second
formula \thetag{10.9} is used for $A\neq 0$ when $\omega_1$ is
defined by \thetag{7.14}.\par
     Let's consider the covariant differential $\theta=\nabla\Theta$.
This is pseudocovectorial field of the weight $-2$ with the following 
components:
$$
\xalignat 2
&\quad\theta_1=\Theta_{\sssize 1.0}-2\,\varphi_1\,\Theta,
&&\theta_2=\Theta_{\sssize 0.1}-2\,\varphi_2\,\Theta.
\hskip-2em
\tag10.10
\endxalignat
$$
Formulas \thetag{10.9} and \thetag{10.10} define the field
$\theta$ effectively in arbitrary coordinates. The quantities
$\varphi_1$ and $\varphi_2$ in \thetag{10.10} should be
calculated either by \thetag{4.8} or \thetag{4.14}. In special
coordinate we can calculate $\theta_1$ and $\theta_2$ explicitly:
$$
\align
&\aligned
 \theta_1=S_{\sssize 1.1}-\frac{6}{5}\,&R_{\sssize 0.2}
 +\frac{6}{5}\,R\,S_{\sssize 1.0}+\frac{12}{5}
 \,Q\,S_{\sssize 0.1}+\\
 \vspace{1ex}
 &+\frac{144}{25}\,S\,Q\,R-\frac{144}{25}\,R
 \,R_{\sssize 0.1}-\frac{324}{125}\,R^3.
\endaligned
\tag10.11\\
\vspace{2ex}
&\theta_2=\frac{486}{125}\,Q\,R^2-\frac{9}{5}\,Q\,S_{\sssize 1.0}
-\frac{108}{25}\,S\,Q^2+\frac{54}{25}\,Q\,R_{\sssize 0.1}-1.
\hskip-4em
\tag10.12
\endalign
$$
Let's raise indices in \thetag{10.11} and \thetag{10.12} by
means of the matrix \thetag{2.5} and let's introduce the
following notations:
$$
\xalignat 2
&\theta^1=\theta_2=C,
&&\theta^2=-\theta_1=D.
\tag10.13
\endxalignat
$$
The quantities $\theta^1$ and $\theta^2$ from \thetag{10.13} 
are components of the pseudovectorial field of the weight 
$-1$. Let's calculate its contraction with $\alpha$:
$$
L=-\frac{5}{9}\,\sum^2_{i=1}\alpha_i\,\theta^i=
\frac{5}{9}\,\sum^2_{i=1}\theta_i\,\alpha^i.
\tag10.14
$$
Pseudoscalar field $L$ from \thetag{10.14} has the weight $0$,
it is connected with the field $K$ by the following relationship:
$$
L=K+\frac{5}{9}.
\tag10.15
$$
Due to \thetag{10.15} and \thetag{10.6} in the fourth case
of intermediate degeneration we have $L\neq 0$. For the
field $\theta$ this means $\theta\nparallel\alpha$. Therefore
we can consider the relationships
$$
\xalignat 2
&\quad\nabla_{\alpha}\alpha=\Gamma^1_{11}\,\alpha+
\Gamma^2_{11}\,\theta,
&&\nabla_{\alpha}\theta=\Gamma^1_{12}\,\alpha+
\Gamma^2_{12}\,\theta,\\
\vspace{-1.7ex}
&&&\tag10.16\\
\vspace{-1.7ex}
&\quad\nabla_{\theta}\alpha=\Gamma^1_{21}\,\alpha+
\Gamma^2_{21}\,\theta,
&&\nabla_{\theta}\theta=\Gamma^1_{22}\,\alpha+
\Gamma^2_{22}\,\theta,
\endxalignat
$$
which are similar to \thetag{5.2}, \thetag{8.3}, and \thetag{9.8}.
For seven coefficients in the expansions \thetag{10.16} we get
the following relationships:
$$
\xalignat 3
&\Gamma^2_{11}=\Gamma^1_{12}=\Gamma^1_{21}=\Gamma^2_{22}=0,
&&\Gamma^1_{11}=\Gamma^2_{21}=-\frac{3}{5}\,N,
&&\Gamma^2_{12}=\frac{12}{5}\,N.
\endxalignat
$$
The formula for eighth coefficients $\Gamma^1_{22}$ is quite 
different:
$$
\aligned
\Gamma^1_{22}&=\frac{9}{25}\,\Theta_{\sssize 0.1}\,Q\,R\,\Theta
+\frac{12}{5}\,Q\,\Theta_{\sssize 0.1}^2
-3\,\Theta_{\sssize 0.1}\,R
-\frac{54}{125}\,R^2\,\Theta^2\,Q-\\
\vspace{1ex}
&-\frac{54}{25}\,R^2\,\Theta
-\frac{9}{5}\,Q\,\Theta\,\Theta_{\sssize 0.2}
-\frac{54}{25}\,Q\,\Theta^2\,R_{\sssize 0.1}
-\Theta_{\sssize 0.2}-\\
\vspace{1ex}
&-\frac{6}{5}\,R_{\sssize 0.1}\,\Theta
+\frac{81}{25}\,S\,Q^2\,\Theta^2+\frac{18}{5}\,S\,Q\,\Theta+S.
\endaligned
\tag10.17
$$
Formula \thetag{10.17} defines pseudoscalar field of the
weight $-4$ in special coordinates. To make effective this 
formula let's use the notations \thetag{10.13} for the
components of the field $\theta$ calculated by \thetag{10.10}:
$$
\aligned
\Gamma^1_{22}=&-\frac{5\,D\,C\,(C_{\sssize 1.0}
-D_{\sssize 0.1})}{9\,L}-\frac{5\,D^2\,C_{\sssize 0.1}
-5\,C^2\,D_{\sssize 1.0}}{9\,L}-\\
\vspace{1ex}
&\qquad-\frac{5\,P\,C^3+15\,Q\,C^2\,D+15\,R\,C\,D^2+5\,S\,D^3}
{9\,L}.
\endaligned
\hskip-3em
\tag10.18
$$\par
     Let's substitute the values of $P$, $Q$, $R$, and $S$
from \thetag{10.2} and \thetag{10.3} into the above formulas 
\thetag{10.8} and \thetag{10.18} for pseudoscalar fields $\Theta$ 
and $\Gamma^1_{22}$. Then
$$
\Theta=\frac{x\,|x|^{-3/4}}{4}\,\sigma(y)+
\frac{4\,x}{3}.
\tag10.19
$$
The expression for $\Gamma^1_{22}$ is more huge:
$$
\aligned
\Gamma^1_{22}=&\frac{|x|^{-3/2}}{64}\,\left(3\,\sigma''(y)\,\sigma(y)
-4\,\sigma'(y)^2\right)+\\
\vspace{1ex}
&+\frac{9\,|x|^{-1/4}}{256}\,\sigma(y)^3+
\frac{3\,x^2\,|x|^{-3/2}}{64}\,\sigma(y)^2.
\endaligned
\tag10.20
$$
Now let's return to the differential equation \thetag{10.7}, which 
we are to make effective. It can be written in the form that doesn't
contain third order derivatives:
$$
\frac{\left(3\,\sigma''(y)\,\sigma(y)-4\,\sigma'(y)^2\right)^3}
{\sigma(y)^{10}}=\const.
\tag10.21
$$
By means of \thetag{10.19} and \thetag{10.20} we compose one more 
pseudoscalar field:
$$
E=\Gamma^1_{22}+\frac{27\,N}{5}\left(\Theta
+\frac{5}{9\,N}\right)^{\! 3}-\frac{3}{4}\left(\Theta
+\frac{5}{9\,N}\right)^{\! 2}.
\tag10.22
$$
Field $E$ in \thetag{10.22} has the weight $-4$. We use it
to construct first scalar invariant for the fifth case of
intermediate degeneration:
$$
I_1=\frac{E^6}{N^8}\left(\Theta+\frac{5}{9\,N}\right)^{\!-20}=
\frac{E^6\,N^{12}}{L^{20}}.
\tag10.23
$$
Using \thetag{10.19} and \thetag{10.20} we can calculate this
invariant \thetag{10.23} in explicit form:
$$
I_1=\frac{6879707136}{390625}\cdot\frac{\left(3\,\sigma''(y)\,
\sigma(y)-4\,\sigma'(y)^2\right)^6}{\sigma(y)^{20}}.
\tag10.24
$$
Comparing formulas \thetag{10.24} and \thetag{10.21}, we can
reformulate theorem~10.1 in the form that admit effective
testing in arbitrary coordinates.
\proclaim{Theorem 10.2} In the fourth case of intermediate 
degeneration algebra of point symmetries of the equation 
\thetag{1.1} is one-dimensional if and only if the invariant
\thetag{10.23} is identically constant. Otherwise this algebra
is trivial.
\endproclaim
\head
11. Fifth case of intermediate degeneration.
\endhead
     In the fifth case of intermediate degeneration we have
the condition $F=0$, which is common for all cases with
degeneration. Parameters $A$ and $B$ do not vanish simultaneously. 
As in the fourth case, here we have the conditions \thetag{6.1} 
and \thetag{10.1}. But the condition \thetag{10.4} is replaced
by quite opposite condition $\sigma(y)=0$ (see \cite{24}). It
can be easily written in effective form:
$$
K+\frac{5}{9}=0.
\tag11.1
$$
Here the field $K$ of the weight $0$ is calculated either by
\thetag{7.22} or \thetag{7.23}. According to the results of
paper \cite{24}, if all above conditions, including \thetag{11.1},
are fulfilled, then the equation \thetag{1.1} can be brought 
to the form
$$
y''=-\frac{5}{12\,x}\,y'+\frac{4}{3}\,x^2\,{y'}^3
\tag11.2
$$
in some special coordinates. This equation \thetag{11.2} 
doesn't contain arbitrary parameters. Therefore its algebra
of point symmetries is quite definite. According to \cite{24},
it is three-dimensional and isomorphic to the matrix algebra
$\sll(2,\Bbb R)$.
\head
12. Sixth case of intermediate degeneration.
\endhead
     Second, third, fourth and fifth cases of intermediate
degeneration are united by the fact that in all these cases
we have $M=0$ and $N\neq 0$ (see conditions \thetag{6.1}).
Sixth an seventh cases are separate in the sense of these
conditions. Here they are replaced by the following one:
$$
N=0.
\tag12.1
$$
The relationship $M=0$ is now derived from \thetag{12.1} and
$F=0$. Parameters $A$ and $B$ do not vanish simultaneously in
sixth and seventh cases too.\par
     Suppose that all above conditions are fulfilled. Let's
transform the equation \thetag{1.1} to the special coordinates 
defined by the conditions \thetag{1.8}. Then, according to the
results of \cite{24}, for $P$ and $Q$ we get
$$
\xalignat 2
&P=0,
&&Q=0.
\tag12.2
\endxalignat
$$
Coefficient $R$ in such coordinates is given by the relationship
$$
R=c(y)\,x+r(y).
\tag12.3
$$
In special coordinates the sixth case of intermediate degeneration 
is distinguished by the additional condition written in terms of
the function $c(y)$ from \thetag{12.3}:
$$
c(y)\neq 0.
\tag12.4
$$\par
     The condition \thetag{12.4}, which is written in special 
coordinates, should be made effective by transformation to the 
arbitrary coordinates. For this purpose we shall use the fact 
that from $F=0$ and $N=0$ we have $M=0$. This means that in sixth 
and seventh cases of intermediate degeneration the pseudoscalar
fields $\Omega$ and $\Lambda$ are defined, as well as the
fields constructed by means of curvature tensor in section 7.
Let's evaluate them in special coordinates:
$$
\xalignat 2
&\Omega=c(y),
&&\Lambda=-2\,c(y).
\tag12.5
\endxalignat
$$
Formulas \thetag{12.5} give us the required effectivization
for the condition \thetag{12.4}:
$$
\Omega\neq 0.
\tag12.6
$$
From \thetag{12.5} and \thetag{12.6} we get equivalent condition
$\Lambda\neq 0$. Moreover, in sixth case of intermediate degeneration
we have the following relationship binding $\Omega$ and $\Lambda$:
$$
\Lambda=-2\,\Omega.
\tag12.7
$$
Let's substitute $Q=0$ from \thetag{12.2} and $\Lambda=-2\,\Omega$ 
from \thetag{12.7} into the formula \thetag{7.8} for the components
of the field $\omega$ in special coordinates. This gives us
$$
\xalignat 2
&\quad\omega_1=\frac{9}{5}\,\Omega,
&&\omega_2=S_{\sssize 1.0}-\frac{6}{5}\,R_{\sssize 0.1}
-\frac{54}{25}\,R^2.
\hskip-3em
\tag12.8
\endxalignat
$$
Let's raise indices and let's introduce the following
notations:
$$
\xalignat 2
&C=\omega^1=\omega_2,
&&D=\omega^2=-\omega_1.
\tag12.9
\endxalignat
$$
From \thetag{12.8} and \thetag{12.6} we see that $\omega$ and $\alpha$
are non-collinear. Therefore we can consider the expansions, which
coincides with \thetag{9.8}:
$$
\xalignat 2
&\nabla_{\alpha}\alpha=\Gamma^1_{11}\,\alpha+
\Gamma^2_{11}\,\omega,
&&\nabla_{\alpha}\omega=\Gamma^1_{12}\,\alpha+
\Gamma^2_{12}\,\omega,\\
\vspace{-1.7ex}
&&&\tag12.10\\
\vspace{-1.7ex}
&\nabla_{\omega}\alpha=\Gamma^1_{21}\,\alpha+
\Gamma^2_{21}\,\omega,
&&\nabla_{\omega}\omega=\Gamma^1_{22}\,\alpha+
\Gamma^2_{22}\,\omega,
\endxalignat
$$
But for the coefficients of the expansions \thetag{12.10} in the
sixth case of intermediate degeneration we get the formulas, which
are different from that of the third case:
$$
\xalignat 2
&\Gamma^1_{11}=\Gamma^2_{11}=\Gamma^2_{12}=0,
&&\Gamma^1_{21}=\Gamma^2_{21}=0,\\
\vspace{1ex}
&\Gamma^2_{12}=1-\frac{12}{25}\,K,
&&\Gamma^2_{22}=\frac{27}{25}\,K.
\endxalignat
$$
The formula for $\Gamma^1_{22}$ is written in special coordinates:
$$
\aligned
\Gamma^1_{22}&=-\frac{93}{25}\,S_{\sssize 1.0}\,R\,\Omega
+\frac{13}{5}\,S_{\sssize 1.0}\,\Omega_{\sssize 0.1}
+S_{\sssize 1.0}
+\frac{306}{25}\,R_{\sssize 0.1}\,R\,\Omega-\\
\vspace{1ex}
&-\frac{78}{25}\,R_{\sssize 0.1}\,\Omega_{\sssize 0.1}
-\frac{6}{5}\,R_{\sssize 0.1}
+\frac{5022}{625}\,R^3\,\Omega
-\frac{702}{125}\,R^2\,\Omega_{\sssize 0.1}-\\
\vspace{1ex}
&-\frac{54}{25}\,R^2
-\frac{9}{5}\,\Omega\,S_{\sssize 1.1}
+\frac{54}{25}\,\Omega\,R_{\sssize 0.2}
+\frac {81}{25}\,\Omega^2\,S.
\endaligned
\hskip-2em
\tag12.11
$$
In arbitrary coordinates we should replace \thetag{12.11} by  
$$
\aligned
\Gamma^1_{22}&=-\frac{5\,C\,D\,(C_{\sssize 1.0}
-D_{\sssize 0.1})}{9\,\Omega}-\frac{5\,D^2\,C_{\sssize 0.1}
-5\,C^2\,D_{\sssize 1.0}}{9\,\Omega}-\\
\vspace{1ex}
&-\frac{5\,P\,C^3+15\,Q\,C^2\,D+15\,R\,C\,D^2+5\,S\,D^3}
{9\,\Omega}.
\endaligned
\tag12.12
$$
Parameters $C$ and $D$ in \thetag{12.12} introduced by
\thetag{12.9} should be calculated by \thetag{7.12} and 
\thetag{7.13} for $B\neq 0$, or by \thetag{7.14} and 
\thetag{7.15} for $A\neq 0$.\par
     In \cite{24} it was shown that at the expense of further
specialization of the choice of coordinates in sixth case of
intermediate degeneration coefficients $R$ and $S$ in the equation
\thetag{1.1} can be brought to the form
$$
\xalignat 2
&\quad R=x,
&&S=x^3+\frac{1}{2}\,x^2+\sigma(y)\,x+s(y).
\hskip-3em
\tag12.13
\endxalignat
$$
There also the following theorem was proved.
\proclaim{Theorem 12.1} In the sixth case of intermediate 
degeneration algebra of point symmetries of the equation 
\thetag{1.1} is one-dimensional if and only if the functions
$\sigma(y)$ and $s(y)$ in \thetag{12.13} are identically 
constant:
$$
\xalignat 2
&s'(y)=0,&&\sigma'(y)=0.
\tag12.14
\endxalignat
$$
If at least one of the conditions \thetag{12.14} fails, then
corresponding algebra of point symmetries is trivial.
\endproclaim
     Let's use \thetag{12.2} and \thetag{12.13} to evaluate
$\Omega$ and $\Lambda$. Their values in special coordinates
appear to be constant:
$$
\xalignat 2
&\Omega=1,
&&\Lambda=-2.
\tag12.15
\endxalignat
$$
The relationships \thetag{12.15} can be derived from \thetag{12.5}
and \thetag{12.13}. Let's also evaluate the field $K$ of the weight 
$0$ and its covariant derivative along $\omega$:
$$
\xalignat 2
&K=x,
&&\nabla_\omega K=\frac{21}{25}\,x^2+x+\sigma(y).
\tag12.16
\endxalignat
$$
Using \thetag{12.16}, we construct another field $L$ of the weight $0$ 
by means of $K$ and $\nabla_\omega K$. We shall choose this field
as the first scalar invariant:
$$
I_1=L=\nabla_\omega K-\frac{21}{25}\,K-K.
\tag12.17
$$
The field \thetag{12.17} is constructed so that its value
in special coordinates coincides with the value of the
function $\sigma(y)$ in \thetag{12.13}:
$$
I_1=L=\sigma(y).
\tag12.18
$$
It's not difficult to evaluate the covariant derivative of $L$
along the vector-field $\omega$:
$$
\nabla_\omega L=-\frac{9}{5}\,\sigma'(y).
\tag12.19
$$
Now we shall use \thetag{12.18}, \thetag{12.19}, and pseudoscalar
field $\Gamma^1_{22}$ of the weight $-2$ from \thetag{12.10} in
order to construct the second scalar invariant:
$$
I_2=\Omega^2\,\Gamma^1_{22}-\nabla_\omega L-\frac{72}{625}\,K^3
+\frac{63}{50}\,K^2+\frac{12}{25}\,K\,L-K-L.
\hskip-3em
\tag12.20
$$
The invariant \thetag{12.20} is constructed so that its value
in special coordinates is proportional to the function $s(y)$ 
in \thetag{12.3}:
$$
I_2=\frac{81}{25}\,s(y).
\tag12.21
$$
From \thetag{12.18} and \thetag{12.21} we derive the following
effectivization for the theorem~12.1 
\proclaim{Theorem 12.2} In the sixth case of intermediate 
degeneration algebra of point symmetries of the equation 
\thetag{1.1} is one-dimensional if and only if both invariants 
$I_1$ and $I_2$ from \thetag{12.17} and \thetag{12.20} 
are identically constant. Otherwise this algebra of point
symmetries is trivial.
\endproclaim
\head
13. Seventh case of intermediate degeneration.
\endhead
     In seventh case of intermediate degeneration the conditions
$F=0$ and $N=0$ remains the same as in sixth case. They give
$M=0$ as immediate consequence. Parameters $A$ and $B$ do not
vanish simultaneously. However, the condition $c(y)\neq 0$ from
\thetag{12.4} is replaced by opposite one: $c(y)=0$. In effective
form this condition is written as the condition opposite to
\thetag{12.6}:
$$
\Omega=0.
\tag13.1
$$\par
     From \thetag{13.1} and \thetag{12.7} we get $\Lambda=0$.
From \thetag{12.8} we find that in special coordinates $\omega_1=0$.
This means that the fields $\omega$ and $\alpha$ are collinear, 
hence the expansions like \thetag{12.10} in this case are impossible.
But at the same time  the collinearity of $\omega$ and $\alpha$ is
exactly the condition for existence of the field $\Theta$ in 
\thetag{10.9}. In special coordinates this field is calculated by
\thetag{10.8}. Repeating the constructions from the section 10,
we define pseudocovectorial field $\theta$ of the weight $-2$
with components \thetag{10.10}. Thereafter we raise indices in
accordance with \thetag{10.13} and we obtain the pseudovectorial
field $\theta$ of the weight $-1$. Formulas \thetag{10.11} and 
\thetag{10.12} for the components of $\theta$ in special 
coordinates in seventh case of intermediate degeneration are 
substantially more simple:
$$
\aligned
\theta^1&=S_{\sssize 1.1}-\frac{6}{5}\,R_{\sssize 0.2}
-\frac{144}{25}\,R\,R_{\sssize 0.1}+\frac{6}{5}\,R
\,S_{\sssize 1.0}-\frac{324}{125}\,R^3,\\
\theta^2&=-1.
\endaligned
\tag13.2
$$
From \thetag{13.2} we see, that $\theta\nparallel\alpha$.
Therefore we can consider the expansions
$$
\xalignat 2
&\quad\nabla_{\alpha}\alpha=\Gamma^1_{11}\,\alpha+
\Gamma^2_{11}\,\theta,
&&\nabla_{\alpha}\theta=\Gamma^1_{12}\,\alpha+
\Gamma^2_{12}\,\theta,\\
\vspace{-1.7ex}
&&&\tag13.3\\
\vspace{-1.7ex}
&\quad\nabla_{\theta}\alpha=\Gamma^1_{21}\,\alpha+
\Gamma^2_{21}\,\theta,
&&\nabla_{\theta}\theta=\Gamma^1_{22}\,\alpha+
\Gamma^2_{22}\,\theta,
\endxalignat
$$
which reproduce \thetag{10.16} for the seventh case of intermediate
degeneration. Almost all coefficients in \thetag{13.3} appear
to be zero. Exception is the coefficient $\Gamma^1_{22}$, which
defines pseudoscalar field of the weight $-4$. In special
coordinates it is calculated by 
$$
\Gamma^1_{22}=S-\Theta_{\sssize 0.2}-\frac{6}{5}\,R_{\sssize 0.1}
\,\Theta-3\,R\,\Theta_{\sssize 0.1}-\frac{54}{25}\,R^2\,\Theta.
\tag13.4
$$
In order to recalculate \thetag{13.4} in arbitrary coordinates
let's denote by $C$ and $D$ the components of the field $\theta$
as it was done in \thetag{10.13}. They should be calculated by
\thetag{10.10}. Then for $\Gamma^1_{22}$ we get
$$
\aligned
\Gamma^1_{22}&=C\,D\,(D_{\sssize 0.1}-C_{\sssize 1.0})
+C^2\,D_{\sssize 1.0}-D^2\,C_{\sssize 0.1}-\\
\vspace{1ex}
&-P\,C^3-3\,Q\,C^2\,D-3\,R\,C\,D^2-S\,D^3.
\endaligned
\tag13.5
$$
Formula \thetag{13.5} has no denominator in right hand side,
since $\Omega=-\theta^1\,\alpha_1-\theta^2\,\alpha_2=1$.\par
     Now let's use the results of the paper \cite{24}. There
it was shown that in the seventh case of intermediate degeneration
upon special choice of coordinates the coefficients of the equation
\thetag{1.1} can be brought to the form
$$
\xalignat 2
&P=0,&&Q=0,\\
\vspace{-1.6ex}
&&&\tag13.6\\
\vspace{-1.8ex}
&R=0,&&S=\frac{1}{2}\,x^2+s(y).
\endxalignat
$$
The conditions \thetag{1.8} for these coordinates appear to be
fulfilled. The structure of algebra of point symmetries for
this case is determined by the function $s(y)$ in \thetag{13.6}.
Here we have two theorems from \cite{24}.
\proclaim{Theorem 13.1} In the seventh case of intermediate
degeneration the algebra of point transformations of the
equation \thetag{1.1} is two-dimensional if and only if
the function $s(y)$ in \thetag{13.6} is identically zero.
This algebra is integrable but it is not Abelian.
\endproclaim
\proclaim{Theorem 13.2} In the seventh case of intermediate
degeneration the algebra of point transformations of the
equation \thetag{1.1} is one-dimensional if and only if
the function $s(y)$ in \thetag{13.6} is nonzero and if
this function is the solution of the following differential 
equation:
$$
4\,s''(y)-\frac{5\,s'(y)^2}{s(y)}=0.
\tag13.7
$$
If $s(y)$ doesn't satisfy the differential equation \thetag{13.7}, 
then the algebra of point symmetries is trivial.
\endproclaim
     Let's substitute \thetag{13.6} into the formula \thetag{10.8} 
for the field $\Theta$, into the formulas \thetag{13.2} for
$\theta$, and into the formula \thetag{13.4} for the field
$\Gamma^1_{22}$. Then we get
$$
\xalignat 2
&\Theta=x,&&\Gamma^1_{22}=\frac{1}{2}x^2+s(y).
\tag13.8
\endxalignat
$$
The components of the field $\theta$ appear to be unitary:
$$
\xalignat 2
&\theta^1=0,&&\theta^2=-1.
\tag13.9
\endxalignat
$$
The field $\Gamma^1_{22}$ in \thetag{13.8} has the weight $-4$,
and the field $\Theta$ has the weight $-2$. Therefore the
following field $L$ has the weight $-4$:
$$
L=\Gamma^1_{22}-\frac{1}{2}\,\Theta^2.
\tag13.10
$$
It is easy to evaluate the field \thetag{13.10} in special 
coordinates:
$$
L=s(y).
\tag13.11
$$
Formula \thetag{13.11} yields the effectivization for the above
theorem~13.1.
\proclaim{Theorem 13.3} In the seventh case of intermediate
degeneration the algebra of point transformations of the
equation \thetag{1.1} is two-dimensional if and only if
the pseudoscalar field $L$ in \thetag{13.10} is identically 
zero. This algebra is integrable but it is not Abelian.
\endproclaim
     In order to make effective the theorem~13.2, first, we note
that for $s(y)\neq 0$ the equation \thetag{13.7} can be written as
$$
\frac{s'(y)^4}{s(y)^5}=\const.
\tag13.12
$$
Then we calculate the covariant derivative of the field $L$ along 
the field $\theta$. The resulting field $\nabla_\theta L$ has the
weight $-5$. Due to \thetag{13.9} its value is given by
$$
\nabla_\theta L=-s'(y).
\tag13.13
$$
Now we are able to compose the scalar invariant $I_1$ by
$\nabla_\theta L$ and $L$:
$$
I_1=\frac{(\nabla_\theta L)^4}{L^5}.
\tag13.14
$$
Weight of tie field \thetag{13.14} is equal to zero. Due to
\thetag{13.11} and \thetag{13.13} the value of this field 
coincides with the left hand side of \thetag{13.12}.
\proclaim{Theorem 13.4} In the seventh case of intermediate
degeneration the algebra of point transformations of the
equation \thetag{1.1} is one-dimensional if and only if
the field $L$ in \thetag{13.10} is nonzero and if the
invariant $I_1$ in \thetag{13.14} is identically constant.
\endproclaim
     If the first hypothesis in the theorem~13.4 is not
fulfilled, i\. e\. if $L=0$, then the algebra of point symmetries
of the equation \thetag{1.1} is two-dimensional due to the 
theorem~13.3. If $L=0$, but $I_1\neq\const$, then this algebra
is trivial due to the theorem~13.2.
\head
14. Acknowledgments.
\endhead
     Author is grateful to Professors E.G.~Neufeld, V.V.~Sokolov,
A.V.~Bocharov, V.E.~Adler, N.~Kamran, V.S.~Dryuma, A.B.~Sukhov 
and M.V.~Pavlov for the information, for worth advises, and for 
help in finding many references below.\par


\Refs
\ref\no 1\by R.~Liouville\jour Jour. de l'Ecole Politechnique,
\vol 59\pages 7--88\yr 1889
\endref

\ref\no 2\by M.A.~Tresse\book Determination des invariants
ponctuels de l'equation differentielle du second ordre
$y''=w(x,y,y')$\publ Hirzel \publaddr Leiptzig\yr 1896
\endref

\ref\no 3\by E.~Cartan\paper Sur les varietes a connection
projective\jour Bulletin de Soc. Math. de France,\vol 52
\yr 1924\pages 205-241
\endref

\ref\no 4\by E.~Cartan\paper Sur les varietes a connexion affine
et la theorie de la relativite generalisee\jour Ann. de l'Ecole
Normale,\vol 40\pages 325--412\yr 1923\moreref\vol 41\yr 1924
\pages 1--25\moreref\yr 1925\vol 42\pages 17-88
\endref

\ref\no 5\by E.~Cartan\paper Sur les espaces a connexion
conforme\jour Ann. Soc. Math. Pologne,\vol 2\yr1923
\pages 171--221
\endref

\ref\no 6\by E.~Cartan\book Spaces of affine, projective and
conformal connection\publ Publication of Kazan University
\publaddr Kazan\yr1962
\endref

\ref\no 7\by G.~Bol\paper Uber topologishe Invarianten von zwei 
Kurvenscharen in Raum\jour Abhandlungen Math. Sem. Univ. Hamburg,
\vol 9\yr 1932\issue 1\pages 15--47
\endref

\ref\no 8\by V.I.~Arnold\book Advanced chapters of the theory
of differential equations\bookinfo Chapter 1, \S~6
\publ Nauka\publaddr Moscow\yr 1978
\endref

\ref\no 9\by N.~Kamran, K.G.~Lamb, W.F.~Shadwick\paper The local
equivalence problem for the equation $d^2y/dx^2\allowmathbreak 
=F(x,y,dy/dx)$ and the Painleve transcendents\jour Journ. of 
Diff. Geometry\vol 22\yr 1985 \pages 139-150
\endref

\ref\no 10\by V.S.~Dryuma\book Geometrical theory of nonlinear
dynamical system \publ Preprint of Math. Inst. of Moldova 
\publaddr Kishinev\yr 1986
\endref

\ref\no 11\by V.S.~Dryuma\paper On the theory of submanifolds
of projective spaces given by the differential equations
\inbook Sbornik statey\publ Math. Inst. of Moldova
\publaddr Kishinev\yr 1989\pages 75--87
\endref

\ref\no 12\by Yu.R.~Romanovsky\paper Calculation of local
symmetries of second order ordinary differential equations
by means of Cartan's method of equivalence\jour Manuscript
\pages 1--20
\endref

\ref\no 13\by L.~Hsu, N.~Kamran\paper Classification of 
ordinary differential equations\jour Proc. of London Math. 
Soc.,\vol 58\yr1989\pages 387--416
\endref

\ref\no 14\by C.~Grisson, G.~Tompson, G.~Wilkens\jour
J.~Differential Equations,\vol 77\yr 1989\pages 1--15
\endref

\ref\no 15\by N.~Kamran, P.~Olver\paper Equivalence problems 
for first order Lagrangians on the line\jour J.~Differential
Equations,\vol 80\yr 1989\pages 32--78
\endref

\ref\no 16\by N.~Kamran, P.~Olver\paper Equivalence of 
differential operators\jour SIAM J.~Math.~Anal.,\vol 20
\yr 1989\pages 1172--1185
\endref

\ref\no 17\by F.M.~Mahomed\paper Lie algebras associated 
with scalar second order ordinary differential equations 
\jour J.~Math. Phys.,\vol 12\pages 2770--2777
\endref

\ref\no 18\by N.~Kamran, P.~Olver\paper Lie algebras of 
differential operators and Lie-algebraic potentials \jour 
J.~Math. Anal. Appl.,\vol 145\yr 1990\pages 342--356\endref

\ref\no 19\by N.~Kamran, P.~Olver\paper Equivalence of 
higher order Lagrangians. I. Formulation and reduction
\jour J.~Math.~Pures et Appliquees,\vol 70\yr 1991
\pages 369--391\endref

\ref\no 20\by N.~Kamran, P.~Olver\paper Equivalence of 
higher order Lagrangians. III. New invariant differential 
equations.\jour Nonlinearity,\vol 5\yr 1992\pages 601--621
\endref

\ref\no 21\by A.V.~Bocharov, V.V.~Sokolov, S.I.~Svinolupov
\book On some equivalence problems for differential equations
\publ Preprint ESI-54, International Erwin Sr\"odinger Institute
for Mathematical Physics\publaddr Wien, Austria\yr1993
\page 12\endref

\ref\no 22\by V.S.~Dryuma\paper Geometrical properties of
multidimensional nonlinear differential equations and phase
space of dynamical systems with Finslerian metric
\jour Theor. and Math. Phys.,\vol 99\issue 2\yr 1994
\pages 241-249
\endref

\ref\no 23\by N.H.~Ibragimov, F.M.~Mahomed\paper Ordinary Differential
Equations\inbook CRC Handbook of Lie Group Analysis of Differential
Equations, Vol. 3\ed N.H.~Ibragimov\publ CRC Press\publaddr
Boca Raton\yr 1996\pages 191-213
\endref

\ref\no 24\by R.A.~Sharipov\paper On the point transformations
for the equations $y''=P+3\,Q\,y'+3\,R\,{y'}^2+S\,{y'}^3$
\jour Electronic Archive at LANL\yr 1997\finalinfo {\bf solv-int
\#9706003}
\endref

\ref\no 25\by M.V.~Pavlov, S.I.~Svinolupov, R.A.~Sharipov\paper
Invariant criterion of hydrodynamical integrability for the
equation of hydrodynamical type\jour Func. analiz i pril. 
\yr 1996 \vol 30\issue 1\pages 18-29
\endref

\ref\no 26\by V.V.~Dmitrieva, R.A.~Sharipov\paper On the point
transformations for the second order differential equations.
\jour Electronic Archive at LANL\yr 1997\finalinfo {\bf solv-int
\#9703003}
\endref

\ref\no 27\by R.A.~Sharipov\book Course of differential geometry 
\publ Publication of Bashkir State University\publaddr Ufa
\yr 1996\page 204
\endref

\ref\no 28\by A.P.~Norden\book Spaces of affine connection
\publ Nauka\publaddr Moscow\yr 1976
\endref
\endRefs
\enddocument
\end